\documentclass[12pt, uplatex]{amsart}
\usepackage{amsmath,amsfonts,amsthm,amscd,amssymb,mathrsfs,amssymb,bm,pb-diagram, here}
\usepackage{setspace}
\usepackage[all,cmtip]{xy}

\usepackage{blkarray}
\usepackage{multirow}
\usepackage{mathtools}
\usepackage{graphicx}

\usepackage{mathrsfs, arydshln}
\usepackage{graphicx}
\usepackage[all]{xy}
\newtheorem{theorem}{Theorem}

\newtheorem{proposition}[theorem]{Proposition}
\newtheorem{lemma}[theorem]{Lemma}
\newtheorem{definition}[theorem]{Definition}

\newtheorem{corollary}[theorem]{Corollary}

\newtheorem{remark}[theorem]{Remark}
\newcommand{\aaa}{\alpha}
\newcommand{\bbb}{\beta}

\newcommand{\id}{{\rm{id}}}
\newcommand{\lmd}{\lambda}
\newcommand{\Lmd}{\Lambda}

\newcommand{\CP}{\mathbb{CP}}

\newcommand{\CC}{\mathbb{C}}

\newcommand{\ZZ}{\mathbb{Z}}
\newcommand{\QQ}{\mathbb{Q}}

\newcommand{\Sing}{{\rm{Sing}\,}}

\newcommand{\qdr}{\CP_1\times\CP_1}

\newcommand{\ol}{\overline}

\newcommand{\lras}{\,\longrightarrow\,}

\newcommand{\Ra}{\Rightarrow}
\newcommand{\set}{\,|\,}

\newcommand{\proofend}{\hfill$\square$}

\newcommand{\inv}{^{-1}}

\newcommand{\Aut}{{\rm{Aut}}}

\newcommand{\Cone}{{\rm{Cone}\,}}

\newcommand{\ms}{\mathscr}
\newcommand{\minus}{\backslash}

\newcommand{\qandq}{\quad{\text{and}}\quad}

\newcommand{\reg}{{\rm{reg}}}

\numberwithin{equation}{section}
\numberwithin{theorem}{section}

\begin{document}
\bibliographystyle{alpha} 
\title[]
{Segre quartic surfaces and minitwistor spaces}
\author{Nobuhiro Honda}
\address{Department of Mathematics, Tokyo Institute
of Technology}
\email{honda@math.titech.ac.jp}

\thanks{The author has been partially supported by JSPS KAKENHI Grant 16H03932.
}
\begin{abstract}
Segre surfaces in the title mean quartic surfaces in $\CP_4$
which are the images of weak del Pezzo surfaces of degree four
under the anti-canonical map.
We first show that minimal minitwistor spaces with genus one 
are exactly Segre quartic surfaces.
By a kind of Penrose correspondence, 
Zariski open subsets of the projective dual varieties of these surfaces
admit Einstein-Weyl structure.
We investigate structures of these dual varieties in detail.
In particular, we determine the degrees of these varieties
(namely the classes of the Segre surfaces),
as well as structure of several components of the boundary divisors
which are the complements of the Einstein-Weyl spaces
in the projective dual varieties.
\end{abstract}

\maketitle
\setcounter{tocdepth}{1}

\section{Introduction}
In a classical terminology, a non-degenerate irreducible 
quartic surface in $\CP_4$ which is not a cone
over a quartic curve in $\CP_3$ nor
a projection of a quartic surface in $\CP_5$
is called a Segre quartic surface.
In a modern language, these are exactly
the images of weak del Pezzo surfaces of degree four
under the anti-canonical maps of the surfaces,
and are realized in $\CP_4$ as complete intersections
of two quadrics.
Segre quartic surfaces are classified into 16 kinds
in terms of normalized quadratic equations of the complete intersections.
Except smooth ones, all these surfaces have isolated singularities,
and all of them are rational double points.

On the other hand, a notion of minitwistor space was originally
introduced by Hitchin in \cite{Hi82-1} as an analogue of 
Penrose's twistor space for a self-dual 4-manifold, and is 
a complex surface which has a smooth rational curve
whose self-intersection number is two.
Such a rational curve is called a minitwistor line,
and they are parameterized by a 3-dimensional complex manifold.
This complex manifold admits a special geometric structure called an Einstein-Weyl structure, 
which is a pair of a conformal structure 
and a compatible affine connection that satisfy 
an Einstein equation.
A smooth quadric and a quadratic cone in $\CP_3$ 
are compact minitwistor spaces,
and minitwistor lines are irreducible hyperplane sections of them.
Essentially these two are all examples of compact minitwistor spaces
in the original sense.

With these backgrounds, 
in the paper \cite{HN11} we showed that, if we allow the
rational curves to have {\em nodes} and at the same time increase
their self-intersection numbers in such a way that 
the parameter space of nodal curves is 3-dimensional, then 
the last space still admits an Einstein-Weyl structure.
We called a complex surface which has such nodal rational curves
a minitwistor space, since if we let the number of nodes zero,
then the surface becomes a minitwistor space in the original sense.
As showed in \cite{HN11},
a compact minitwistor space in the generalized sense
is always a rational surface.
Hence, all minitwistor lines on it are linearly equivalent to each other.
If we write $g$ for the number of nodes of minitwistor lines,
then the linear system generated by minitwistor lines on a compact 
minitwistor space is $(3+g)$-dimensional,
and the map induced by this linear system is
always a birational morphism over the image.
We call a compact minitwistor space to be {\em minimal}\,
if this birational morphism is an isomorphism over the image.
Hence, compact minimal minitwistor spaces are
projective surfaces in $\CP_{3+g}$,
and minitwistor lines on them are obtained as hyperplane
sections of the surfaces.
So the situation is quite similar to the classical case of $g=0$,
but when $g>0$ a generic member of the linear system
is not a minitwistor line because such a member is a smooth
curve whose genus is exactly $g$.
By this reason in this article
we call the number of nodes of minitwistor lines
the {\em genus} of a minitwistor space.

In the present article 
we shall investigate compact minimal minitwistor spaces whose genus are one
and the associated Einstein-Weyl spaces.
Our first main result means that
such minitwistor spaces are exactly the Segre quartic surfaces 
(Theorem \ref{thm:Segmt}).
This implies that, in contrast to the classical case of $g=0$,
there are a variety of compact minitwistor spaces with $g=1$, 
but they can still be concretely classified.
Minitwistor lines on a Segre quartic surface $S\subset\CP_4$ are
hyperplane sections which are tangent to $S$
at exactly one point.
Precisely speaking, the tangency means that the hyperplane contains a tangent plane
of $S$ at some smooth point of $S$,
and the minitwistor line cut out by the hyperplane has a node at the tangent point.
A completion of the space of such hyperplanes
is nothing but the projective dual variety of the Segre surface.
Hence, the Einstein-Weyl spaces associated to the Segre quartic surfaces
are realized as Zariski-open subsets of the dual varieties
of the surfaces.
From a reflexibility for the operation
of taking the projective dual,
a Penrose type correspondence between compact minimal minitwistor spaces
with genus one and the associated Einstein-Weyl spaces
can be understood as a projective duality.
The complements of the Einstein-Weyl spaces in the dual varieties
of Segre surfaces are of special interest, 
because
in the case of the smooth quadric $Q$ in $\CP_3$,
the complement is exactly the dual quadric $Q^*\subset\CP_3^*$,
and this can be regarded as a complexification of the conformal infinity of
the real form of the complex Einstein-Weyl space.
We call 2-dimensional components of the above complements in 
the dual variety
 {\em divisors at infinity}.

A large part of this article is devoted to investigate 
the dual varieties of Segre quartic surfaces and the divisors at infinity on
the dual varieties.
We show that all these dual varieties are rational (Proposition \ref{prop:rat}).
Next in Section \ref{ss:cf}, we give a formula which expresses the degrees
of the dual varieties of any Segre surfaces in terms of 
the types of singularities of the surfaces (Theorem \ref{thm:cf}).
The degree of the dual variety will be smaller as the singularities
of the surface become complicated.
In Section \ref{ss:segs}, 
we recall Segre symbol by which all Segre quartic surfaces can be
classified in a systematic way,
and complete a classification of Segre surfaces (Proposition \ref{prop:cmp}).
In Section \ref{ss:dcS}, we investigate double covering structure on many kinds of Segre surfaces, which can be detected from the Segre symbols.
The base space of the double covering is either a smooth quadric
or a quadratic cone in $\CP_3$
(Propositions \ref{prop:SS} and \ref{prop:SC}).
In Section \ref{ss:sil} we investigate divisors at infinity on the dual varieties
of Segre surfaces.
We will find two kinds of such divisors.
A first kind arises from straight lines lying on Segre surfaces,
and the divisors are 2-planes (Proposition \ref{prop:line}).
The second kind  comes from the double covering 
structure over a smooth quadric given in Section \ref{ss:dcS}, 
and they are smooth quadric surfaces which are dual to the last smooth quadrics
(Proposition \ref{prop:dcQ}).
We determine the numbers of these divisors at infinity,
for each type of Segre quartic surface (Tables \ref{table1},
\ref{table2} and \ref{table3}).
Furthermore, by using deformation theoretic argument,
we show that the dual varieties of Segre surfaces
have ordinary double points along these two kinds of divisors
(Proposition \ref{prop:si1}).
In other words, the dual varieties of Segre quartic surfaces
have self-intersection along the 2-planes or the smooth quadrics.

As above, the double covering structures on the Segre surfaces
are induced by a projection from a point,
and the point determines the dual hyperplane in the dual projective space $\CP_4^*$.
In Section \ref{s:shps}, with the aid of the results
in the previous section, we determine the sections of the dual varieties
of the Segre surfaces by these hyperplanes, in explicit forms
(Propositions \ref{prop:swqb},  \ref{prop:swqb2} and  \ref{prop:swqb3}).
For some Segre surfaces, these descriptions give another proof 
for the class formula from different angle.
Next, in Section \ref{ss:dgn}, we discuss several typical transitions
between Segre quartic surfaces.
Finally, in Section \ref{ss:cr}, we give a remark about 
null surfaces in some of the present Einstein-Weyl spaces,
and also provide a question about present minitwistor spaces,
in connection with twistor spaces associated to self-dual 4-manifolds.

\section{Minitwistor spaces with genus one and Segre quartic surfaces}\label{s:mts}
We begin with the definition of minitwistor spaces
in the sense of \cite{HN11}.
This naturally includes minitwistor spaces in the original sense given in \cite{Hi82-1}
as the simplest case.
Let $g\ge 0$ be an integer.
By a {\em $g$-nodal rational curve}, we mean a
rational curve which has exactly $g$  
nodes (i.e.\,ordinary double points) as its 
all singularities.
In particular, when $g=0$, it is just a smooth rational curve.

\begin{definition}\label{def:mt}
{\em
Let $g\ge 0$ be an integer and $S$ 
a normal, compact and irreducible complex surface.
A $g$-nodal rational curve $C$ 
on $S$ is called a {\em minitwistor line}
if it is contained in $S_{\reg}$,  the smooth locus of $S$,
and the self-intersection number satisfies $C^2 = 2 + 2g$.
A normal compact complex surface having 
a minitwistor line with $g$ nodes as above
is called a {\em minitwistor space with genus $g$}.
}
\end{definition}

By \cite[Proposition 2.6]{HN11}, any minitwistor space in this sense is a rational surface.
Further by \cite[Proposition 2.8]{HN11},
the linear system $|C|$ generated by minitwitor lines
on a minitwistor space $S$ as above
is $(3+g)$-dimensional and base point free.

\begin{proposition}\label{prop:genus}
In the situation of the above definition,
a generic member of the system $|C|$ is smooth and 
it is a curve of genus $g$.
\end{proposition}

\proof
The smoothness of a generic member of $|C|$ follows from Bertini's theorem.
Let $C\subset S$ be a minitwistor line, 
$\nu:\tilde S\to S$ the blowup of $S$ at
the $g$ nodes of $C$,
and $\tilde C$ the strict transform of $C$ into $\tilde S$.
Evidently the curve $\tilde C$ is a smooth rational curve.
If $E_1,\dots, E_g$ are the exceptional curves
of $\nu$, we have
\begin{align}
K_{\tilde S}\cdot\tilde C &= 
\Big(
\nu^*K_S + \sum_{i=1}^g E_i
\Big)\cdot
\Big(
\nu^*C -2 \sum_{i=1}^g E_i
\Big)= K_S\cdot C + 2g
\label{adjn1}
\end{align}
and 
\begin{align}
\tilde C^2 &= 
\Big(
\nu^*C -2 \sum_{i=1}^g E_i
\Big)^2= C^2 - 4g = 2-2g.
\label{adjn2}
\end{align}
Further, since $\tilde C$ is a smooth rational curve,
we have by adjunction 
\begin{align*}
K_{\tilde S}\cdot \tilde C + \tilde C^2 = -2.
\end{align*}
Substituting \eqref{adjn1}
and \eqref{adjn2} to this equality,
we obtain $K_S\cdot  C=-4$.
Hence, if $C'$ is a smooth member of $|C|$, 
we have $K_S\cdot  C'=-4$.
Hence, again by adjunction and using $(C')^2 = C^2 = 2+2g$, we obtain that
the genus of the curve $C'$ is exactly $g$.
\proofend

\medskip
This is why we call $g$ the genus of a minitwistor space.
Thus minitwistor lines are obtained as a degeneration 
of smooth members of the linear system $|C|$ into
irreducible curves which have exactly $g$ nodes.
We denote the morphism associated to the $(3+g)$-dimensional system $|C|$ by
\begin{align}\label{phi1}
\phi:S\lras\CP_{3+g}.
\end{align}
By \cite[Proposition 2.8]{HN11}, this morphism 
is always birational over the image.
However, the morphism $\phi$ can contract curves on $S$.
Indeed, any blow-up of a minitwistor space with genus $g$ is again a minitwistor space with genus $g$,
but the map associated to the system generated by minitwistor lines
contracts the exceptional curve of the blow-up.
In order to ignore these redundant spaces, 
we call that a minitwistor space $S$ is {\em minimal}\,
if the birational morphism $\phi$ in \eqref{phi1}
is an embedding,
so that it does not contract any curves.
Note that this is little stronger than the minimality 
introduced in \cite{HN11} to the effect that 
in \cite{HN11} the minimality means that
the curves contracted by $\phi$ were $(-1)$-curves only,
while in the present paper we are allowing 
other curves (e.g.\,a $(-2)$-curve which
would exist when we further blow up a point on the 
exceptional curve of the first blow up) to be contracted.
The present definition seems more natural
for the purpose of dropping these redundant spaces.
Of course, the minimality does not imply that
$S$ does not have a $(-1)$-curve.

Since the system $|C|$ is base point free
and the associated morphism $\phi$ is birational 
over the image as above, for the degree of the image surface
$\phi(S)\subset \CP_{3+g}$, we always have
\begin{align}\label{deg1}
\deg\phi(S) = C ^2 = 2+ 2g.
\end{align}

When $g=0$, the conditions in Definition
\ref{def:mt} mean that $C$ is a
smooth rational curve satisfying $C^2 = 2$,
and this agrees with the original definition
of a minitwistor space given in \cite{Hi82-1}.
In this case, the image surface $\phi(S)\subset\CP_4$ is either a smooth quadric or the cone over an irreducible conic.
These are all examples of minitwistor spaces classically known, 
and minimal minitwistor spaces with genus zero
are exactly these two surfaces \cite[Proposition 2.14]{HN11}.

Next we define a Severi variety of rational curves
associated to a minitwistor space
in the present sense.

\begin{definition}\label{def:Severi}
{\em
Let $S\subset\CP_{3+g}$ be a minimal minitwistor space with genus $g>0$.
Let $W_0$ be the subset of the dual projective space 
$\CP^*_{3+g}$ consisting of hyperplanes 
$H\subset\CP_{3+g}$ such that the hyperplane sections $S\cap H$ are minitwistor lines 
(with $g$ nodes).
We write $W$ for the closure of $W_0$, taken in $\CP^*_{3+g}$.
This is a subvariety in $\CP^*_{3+g}$, and
we call it the {\em Severi variety of $g$-nodal rational curves} on $S$.
}
\end{definition}

The space $W_0$ is a Zariski-open subset of
the closure $W$.
In other words, the Severi variety $W$ is a compactification of $W_0$.
We note that in the above definition of the space $W_0$,
since we are requiring that a minitwistor line
is included in $S_{\reg}$ as in Definition \ref{def:mt},
we are requiring that the hyperplane sections $S\cap H$
do not pass any singularity of $S$.
If a hyperplane $H$ passes a singularity of $S$,
say $p$, 
then the section $S\cap H$ is always singular
at $p$, and the section $S\cap H$ can be a nodal 
rational curve with a correct number of nodes.
But in general such a section does not admit an equisingular displacement in $|C|$ 
which avoids the singularity of $S$,
and in that case $H\not\in S^*$.
See Remark \ref{rmk:mtl} for concrete examples of this kind.

\begin{remark}{\em
Although it seems quite likely, we do not know 
whether the Zariski open subset $W_0$ in $W$ is 
precisely the smooth locus of the Severi variety $W$.
}
\end{remark}
Now from \cite[Proposition 2.6 and Theorem 2.10]{HN11},
as a natural generalization of a well-known result in the case of genus zero given in \cite{Hi82-1}, we have
\begin{proposition}\label{HNMT}
The subset $W_0\subset\CP_{3+g}^*$ is a 3-dimensional complex manifold, and it admits a torsion-free Einstein-Weyl structure.
\end{proposition}

In this article, we are mainly interested in
minimal minitwistor spaces with genus one.
In this case the Severi variety $W$ is nothing
but an object which is used more often.
To explain this object precisely, let $S$ be an irreducible,
non-denenerate 2-dimensional
subvariety in $\CP_4$.
We write $S^*\subset\CP_4^*$ for the dual variety of $S$.
This is by definition \cite[Def.\,1.1]{Tev}
the closure in the dual projective space
$\CP_4^*$ of the locus of hyperplanes 
which contain tangent planes of $S$
at some smooth points of $S$.
(If $S$ is smooth, we do not need to take the closure.)
Let $I(S)$ be the incidence variety for $S$ and $S^*$.
Namely $I(S)$ is the closure of the set
$$
\big\{ (p,H)\in \CP_4\times\CP_4^*\set
p\in S_{\reg},\,\,T_p S\subset H
\big\},
$$
taken in $\CP_4\times\CP_4^*$.
This is a subvariety in $\CP_4\times\CP_4^*$,
and the dual variety $ S^*$ 
is the image of $I(S)$ under the  
the projection to the second factor $\CP_4^*$.
Hence, there is a double fibration 
\begin{align}\label{diagram:fac1}
 \xymatrix{ 
&I(S) \ar[dl] \ar[dr] &\\
S && S^*
 }
\end{align}
Over the smooth locus $S_{\reg}$ of $S$, the incidence variety $I(S)$ is
a fiber bundle whose fibers are projective lines.
In particular, $I(S)$ is 3-dimensional.
Also, since $S$ is supposed to be irreducible, $I(S)$ is always an irreducible variety.
Therefore, so is the dual variety $S^*$.

Now suppose that $S$ is a minimal minitwistor space
with genus one.
So the morphism \eqref{phi1} provides
a projective embedding $S\subset\CP_4$.
If $H\subset\CP_4$ is a hyperplane 
which belongs to the subset $W_0$, 
the intersection $S\cap H$ is a minitwistor line
from Definition \ref{def:Severi}.
If $p$ is the node of this minitwistor line
then we have $T_pS\subset H$ because otherwise 
$S\cap H$ would be smooth at the point $p$.
Thus $H\in W_0$ implies $H\in S^*$.
Namely $W_0\subset S^*$.
Therefore, since $S^*$ is closed, we obtain
$W \subset S^*$ for the closure $W$ of $W_0$.
Since $S^*$ is irreducible as above and at most 3-dimensional, this means 
that $W=S^*$.
Thus {\em in the case of a minitwistor space with genus one, 
the Severi variety $W$ in Definition \ref{def:Severi} is nothing but the dual variety
$S^*$, and it is always irreducible.}

\begin{remark}\label{rmk:mtl}{\em
From this we can easily obtain an example of a 1-nodal rational curve
on a minitwistor space with genus one, which does not belong to the 
Zariski-open subset $W_0$ (namely which is not a minitwistor line)
as follows.
Take any minimal minitwistor space $S\subset\CP_4$ with genus one which has at least one
ordinary double point $p$.
(See Section \ref{ss:dcS} for examples such minitwistor spaces.)
Let $p^*\subset\CP_4^*$ be the hyperplane which is dual to $p$.
So a hyperplane $H\subset\CP_4$ belongs to $p^*$ iff $p\in H$.
By the irreducibility of $S^*$, we have  $p^*\not\subset S^*$.
Therefore, the hyperplane section $S^*\cap p^*$ is a divisor on $p^*=\CP_3$.
If we take a generic $H\in p^*$ which does not belong to 
$S^*$, the section $S\cap H$ is a 1-nodal rational curve
whose node is exactly $p$,
and it cannot be deformed into a minitwistor line by any small displacement in $S$
since $H\not\in S^*$.
This implies that $S\cap H$ is not a minitwistor line.
}
\end{remark}

As is already mentioned, the structure of minimal minitwistor spaces
with genus zero is strongly constrained.
Next we would like to see that a constraint
for the structure of minimal minitwistor spaces
with genus one is somewhat moderate, 
and they are exactly particular surfaces that are classically known.
In order to explain what are these surfaces
and put them in a broader context,
we discuss classical results on del Pezzo surfaces
of arbitrary degrees.

Suppose $n\ge 3$. All surfaces in $\CP_n$ appearing  below 
are assumed to be irreducible and non-degenerate,
but smoothness and even normality are not assumed.
It is classically known \cite[p.\,174]{GH} that the degree of
any surface $S$ in $\CP_n$ 
is bounded from below as
$$
\deg S \ge n-1.
$$
Surfaces which attain the minimal degree
$(n-1)$ are classically classified \cite[p.\,525]{GH}.
Surfaces with the second smallest degree,
namely surfaces with degree $n$ in $\CP_n$
are also classified, and they belong to one of the following 
three kinds of surfaces \cite[Section 8.1]{Dol}:
\begin{itemize}
\item[(a)] the images of  surfaces
of degree $n$ in $\CP_{n+1}$ under projections
from points,
\item[(b)] the cones over irreducible curves
of degree $n$ in $\CP_{n-1}$,
\item[(c)] surfaces not included in (a) nor (b).
\end{itemize}
A difference between these surfaces
is that  generic hyperplane sections
of surfaces in (a) have arithmetic genus zero, while 
they are one for surfaces in (b) and (c).
See \cite[Section 8.1]{Dol}.
So an irreducible non-degenerate surface
of degree $n$ in $\CP_n$ belongs to
(c) iff a generic hyperplane section has arithmetic genus one and the surface is not 
the cone over a  curve of degree $n$ in $\CP_{n-1}$.

According to \cite[Definition 8.1.5]{Dol}, in a classical terminology, 
an irreducible non-degenerate surface in $\CP_n$
belonging to the case (c) is called a {\em del Pezzo surface} of degree $n$.
The degree of a del Pezzo surface is at most 9 \cite[Propositions 8.1.7 and 8.1.8]{Dol}.
Any del Pezzo surface is normal and has at worst
rational double points.
For smooth ones, these surfaces are exactly del Pezzo surfaces in modern definition.
For singular ones, 
the minimal resolutions of del Pezzo surfaces
in this classical sense have $(-2)$-curves, and in modern language they are often called weak del Pezzo surfaces.

Classically del Pezzo surfaces of degree four
are called {\em Segre quartic surfaces}
(\cite[Section 8.6]{Dol}).
By \cite[Theorem 8.6.2]{Dol}, any Segre quartic surface
is a complete intersection of two quadrics in $\CP_4$.
From this we readily obtain that the anti-canonical class of 
a Segre quartic surface is the class of hyperplane sections.
Any Segre quartic surface has at most finite number of lines on it
(see \cite[Section 8.6.3]{Dol}).
Note that a complete intersection of two quadrics
in $\CP_4$ is not necessarily a Segre surface 
because the cone over a quartic curve in $\CP_3$ 
is typically a complete intersection of two quadrics
in $\CP_4$.

With these preliminaries, we have

\begin{theorem}\label{thm:Segmt}
Any minimal minitwistor space
with genus one is a Segre quartic surface.
Conversely, any Segre quartic surface is 
a minimal minitwistor space with genus one.
\end{theorem}

\proof
Let $S$ be a minimal minitwistor space with genus one.
If $C$ is a minitwistor line on $S$, as is already remarked,
the surface $S$ is embedded in $\CP_4$
by the complete linear system $|C|$,
and in particular $S$ is non-degenerate
in $\CP_4$.
Since $C^2 = 2 + 2g = 4$, we have $\deg S = 4$.
Thus, $S$ is a non-degenerate irreducible 
quartic surface in $\CP_4$.
Hence, by letting $n=4$ in the above classification
of surfaces in $\CP_n$ of degree $n$,
in order to show that $S$ is a Segre surface,
it is enough to show that the surface $S$ does not belong
to the above classes (a) nor (b).
%
The class (a) is immediately rejected
since the arithmetic genus of 
a generic hyperplane section of surfaces in (a) is zero, while it is one for minitwistor spaces with genus one by Proposition \ref{prop:genus}.
Next, we show that the cones as in the class (b) cannot be
a minitwistor space.
Recall that we are supposing normality for a
minitwistor space (see Definition \ref{def:mt}).
Therefore, because the cone over a singular curve
is non-normal, the cone can be a minitwistor space with genus one
only when the base curve is smooth.
But even in that case the cone cannot be 
a minitwistor space with genus one
because a hyperplane section of the cone is singular only
when it passes the vertex of the cone,
and a section of the cone by such a hyperplane
consists of generating lines of the cone and therefore it
 cannot be a minitwistor line.
Hence the surface $S$ belongs to the class (c), and this means that $S$ has to be a Segre quartic surface.

To prove the converse, let $S\subset\CP_4$ be a Segre quartic surface.
We first show that the dual variety $S^*$ of $S$ is 3-dimensional.
Recall that a surface in $\CP^4$ is said to be ruled if any point
of the surface is passed by a line lying on the surface.
It is easy to see that the dual variety of any ruled surface is at most 2-dimensional.
Moreover, this can happen only for ruled surfaces \cite[Theorem 1.18]{Tev}.
But any Segre quartic surface is not ruled because it has only a finite number of lines on it.
Therefore, $\dim S^*=3$.
As in \cite[Section 2.1]{VoII}, that
the dimension of the dual variety is maximal
(i.e. it is a hypersurface in the dual projective space)
 implies that a hyperplane section of $S$
which corresponds to a generic element of the dual variety $S^*$ has
precisely one node as its all singularity.
From genericity such a hyperplane section can be assumed not to pass any singularity of $S$,
and if $C$ is such a hyperplane section, we have $C^2 = H.H.S = 4$.
Moreover, the curve $C$ is rational since it has
exactly one node as its all singularity and
its arithmetic genus is one.
Therefore, $C$ is a minitwistor line with one node, and hence the surface $S$ is 
a minimal minitwistor space with genus one.
\proofend

\medskip
Thus, minimal minitwistor spaces with genus one
are exactly Segre quartic surfaces.
Hence from Proposition \ref{HNMT}, 
the open subset $W_0$ of the dual variety (or equivalently of Severi variety) in Definition \ref{def:Severi}
of any Segre surface
admits a torsion-free Einstein-Weyl structure.

For a birational property of the dual varieties
of Segre quartic surfaces, we have the following

\begin{proposition}\label{prop:rat}
The dual variety of any Segre quartic surface is rational.
\end{proposition}

\proof
Let $S\subset\CP_4$ be a Segre quartic surface.
Recall that, writing $I(S)$ for the incidence variety as before,
we have the diagram
\begin{align}\label{diagram:fac2}
 \xymatrix{ 
&I(S) \ar[dl] \ar[dr] &\\
S && S^*
 }
\end{align}
Take and fix a point $q$ of $\CP_4$ in which $S$ is embedded.
Then for a generic
point $p\in S_{\reg}$, we have $q\not\in T_pS$. 
Hence the linear subspace  $H(p)$ spanned by 
a 2-plane $T_pS$ and the point $q$
is a hyperplane, and it contains $T_pS$.
Hence the pair $(p,H(p))\in \CP_4\times\CP_4^*$ belongs to $I(S)$.
So the assignment $p\mapsto (p,H(p))$ defines a
rational section of the projection
$I(S)\to S$.
Generic fibers of this projection are projective lines,
and therefore, the presence of 
the rational section implies \cite[Lemma 3.4]{CK98}
that $I(S)$ is birational to the product $S\times\CP_1$.
Moreover, as mentioned right after Definition \ref{def:mt}, any minitwistor space is a rational surface.
Hence the variety $I(S)$ is rational.
Since  the variety $S^*$ is the image of $I(S)$ under the 
projection to $\CP_4^*$, this means that $S^*$ is unirational.
For the rationality, it is enough to show that 
the projection $I(S)\to S^*$ is of degree-one.
If the degree is greater than one,
a hyperplane $H\in S^*$ which is generic in $S^*$
contains at least two distinct tangent spaces of $S$.
This means that $H\cap S$ has at least two singularities.
But as mentioned in the last part of the proof of Theorem
\ref{thm:Segmt}, the intersection $H\cap S$ has exactly one node
as its all singularities for generic $H\in S^*$.
Hence the projection $I(S)\to S^*$ is of degree-one.
\proofend

\medskip
Note that the proof works for any rational variety $X\subset\CP_N$ 
which is not ruled.

\section{Study on the dual varieties of Segre quartic surfaces}
\subsection{The class formula for Segre quartic surfaces}
\label{ss:cf}
In this subsection we determine the degrees of
the dual varieties of Segre quartic surfaces, namely the {\em classes}
of the surfaces, by blowing up the surfaces
at the intersection with a generic 2-plane in $\CP_4$
and then investigating singular fibers of the resulting elliptic
fibrations induced on the blowups.

First, we explain how
one can obtain an elliptic fibration.
Let $S\subset\CP_4$ be any Segre quartic surface,
and $P\subset\CP_4$ a 2-plane which is sufficiently general so that 
$P$ intersects $S$ transversally at 
any point of $S\cap P$.
This in particular means that $S\cap P$
consists of four points as  $\deg S=4$, 
and $S$ is smooth at these four points.
If $p$ is any one of these points,
from the transversality
we have $T_pS\cap T_pP = 0$ for the intersection of 
tangent spaces.
Therefore, for any hyperplane $H$ containing $P$,
we have $T_pS \not\subset H$, and 
hence 
the hyperplane section $$S_H:=S\cap H$$ is smooth at
any of the four points $S\cap P$.
These hyperplane sections are quartic curves
in $H=\CP_3$,
and belong to the anti-canonical class of $S$.
Further for any two distinct hyperplanes containing $P$, the anti-canonical curves
intersect transversally at the four points
$S\cap P$.
Let $\CP_4\to \CP_1$ be the generic projection from 
the 2-plane $P$.
Fibers of this projection are hyperplanes that contain
$P$,
and by intersecting with $S$, we obtain a pencil of anti-canonical curves on $S$.
This pencil has the four points $S\cap P$ as the base locus.
Let $S'\to S$ be the blowing-up at the four points
$S\cap P$. Equivalently the surface $S'$ is the strict transform of 
$S$ under the blowing-up of $\CP_4$ along the 2-plane $P$.
By the transversality the base points of the above pencil on $S$
are eliminated through the blowup, and we obtain a morphism.
We write it as 
\begin{align}\label{fprime}
f':S'\to \CP_1.
\end{align}
All fibers of $f'$ are isomorphic to the corresponding members
of the original pencil on $S$.
By Bertini's theorem
fibers of $f'$ are smooth except for finite ones.
Obviously all fibers of $f'$ are 
anti-canonical curves on $S'$.
In particular it is an elliptic curve
as long as it is smooth.
Thus the morphism $f':S'\to \CP_1$ is an elliptic fibration.
Since all blown-up points on $S$ are smooth points of $S$,
the surface $S'$ has the same singularities
as $S$ has.

Let $\tilde S'\to S'$ be the minimal resolution 
of all singularities of $S'$.
If $S$ is smooth we promise $\tilde S'=S'$.
Since all singularities of $S'$ are rational double points,
all components of the exceptional divisors of the resolution are $(-2)$-curves.
We write 
\begin{align}\label{fprimet}
\tilde f':\tilde S'\to \CP_1
\end{align}
for the composition $\tilde S'\to S'
\stackrel{f'}\to\CP_1$.
This is also an elliptic fibration,
but this time $\tilde S'$ is smooth.
Since all singularities of $S'$ are rational double points, fibers of $\tilde f'$ are still anti-canonical curves on $\tilde S'$,
and therefore we have $K^2 = 0$ for the surface $\tilde S'$.
By Hartogs theorem this readily means that any fiber of $\tilde f'$ does not contain a $(-1)$-curve.
Namely the elliptic fibration $\tilde f'$ in \eqref{fprimet} is relatively minimal.

Obviously, these constructions depend only 
on the choice of the 2-plane $P\subset\CP_4$,
where the 2-plane has to satisfy the transversality
for intersection with $S$.
But for the purpose of calculating the degrees of 
the dual varieties of the Segre surfaces,
we need to choose the plane $P$ more carefully to make singular
fibers of the elliptic fibration
\eqref{fprime} or equivalently \eqref{fprimet}
in most generic forms.
For this purpose, we show
%
%

\begin{lemma}\label{lemma:I1}
There exists a plane $P\subset\CP_4$ such that 
the (singular) elliptic fibration
$f':S'\to\CP_1$ 
induced by $P$ as in \eqref{fprime}
satisfies the following properties.
\begin{itemize}
\item[(i)]
Two singularities of the surface $S'$ do not belong to 
the same fiber of $ f'$.
\item[(ii)]
Any singular fiber of $f'$ on which a singularity
of $S'$ belongs has no singularity other than
that singularity of $S'$.
\item[(iii)]
If a singular fiber of $ f'$ has no singularity
of $S'$ on it,
the fiber is of type {\rm I}$_1$.
\end{itemize}
\end{lemma}


\proof
This is a problem of the existence of a line 
in the dual space $\CP_4^*$ which defines
a pencil on $S$ whose associated morphism
$f':S'\to\CP_1$ satisfies the three
genericity conditions in the lemma.

First let $\ms S$ be the set of hyperplanes in $\CP_4$
such that the sections of $S$ by the hyperplanes are singular.
By Bertini's theorem $\ms S$ is a strict subvariety of the dual space $\CP_4^*$.
Obviously the dual variety $S^*$ is an irreducible component of 
the subvariety $\ms S$.
Also if $p_1,\dots, p_k$ are 
all singularities of $S$, 
the dual hyperplanes $p_1^*,\dots, p_k^*\subset\CP_4^*$,
which are the 
sets of hyperplanes in $\CP_4$ which pass the singular points $p_1,\dots, p_k$
respectively, are components of $\ms S$.
To show that these are all components
of $\ms S$,
let $H\subset\CP_4$ be any hyperplane such that $
S_H:=S\cap H$ is singular, and suppose that $H$ does not
belong to $S^*\cup(p_1^*\cup\dots \cup p_k^*)$.
If there is a singularity, say $p$, 
of $S_H$ which is a smooth point of $S$,
we have $T_pS\subset H$ and from the definition of the dual variety $S^*$, this means $H\in S^*$
which contradicts our choice of $H$.
Hence, all singularities of $S_H$ have to be
singular points of $S$.
This implies $H\in p_i^*$ for some $i$,
and again this cannot happen from our choice of $H$.
Hence  we have
\begin{align}\label{Sdp}
\ms S = S^*\cup \big(p_1^*\cup\dots \cup p_k^*\big).
\end{align}

For a portion of the locus in the dual space $\CP_4^*$
which should be avoided from a line to pass,
for any different indices $i,j\le k$,
we put $P_{ij}:=p_i^*\cap p_j^*$.
This is a 2-plane in $\CP_4^*$ and 
is the space of hyperplanes which pass
$p_i$ and $p_j$.
When the surface $S$ has at most one singularity,
we do not need these in the following argument.
In particular, for any $H\in P_{ij}$,
the section $S_H$ has singularities
at least at $p_i$ and $p_j$.
For each index $i\le k$, let $D_i$ be the subset
of the dual hyperplane $p_i^*$ such that if $H\in D_i$,
the section $S_H$ has a singularity other than $p_i$.
We show that for each $i=1,\dots,k$,
\begin{align}\label{olD}
\ol D_i = (S^*\cap p_i^*) \cup \Big(
\bigcup_{j\neq i} P_{ij}
\Big)
\end{align}
holds, where $\ol D_i$ is the closure of $D_i$
in $\CP_4^*$.
 The inclusion `$\supset$' is obvious.
For the reverse inclusion,
take any $H\in D_i$ which does not belong to
$P_{ij}$ for any $j\neq i$.
Then the section $S_H$ has a singularity 
which is necessarily a smooth point of $S$.
This means $H\in S^*$ and hence, $H\in S^*\cap p_i^*$.
Therefore, $D_i$ is included in RHS of \eqref{olD}.
Since $\ol D_i$ is obviously a subvariety of $\CP_4^*$
and in particular closed, this implies the inclusion `$\subset$'
in \eqref{olD}. Thus, the equality \eqref{olD} holds.

Next, as in Definition \ref{def:mt}, let $W_0$ be the subset of $S^*$ 
formed by  hyperplanes $H\subset\CP_4$ such that 
$S\cap H$ is a minitwistor line (with one node
in the present situation).
From the definition of a minitwistor line, 
the node  is not a singularity of $S$.
We denote $A:=S^*\minus W_0$. 
Since $W_0$ is Zariski-open in $W=S^*$ and non-empty,
$A$ is a strict subvariety of $S^*$.
In particular any irreducible component of $A$
is at most 2-dimensional.

Finally, we consider the set of 2-planes in $\CP_4$
which do not intersect transversally with $S$.
This is a strict subvariety of the Grassmanian of 
2-planes in $\CP_4$.
Let $T$ be the complement of this subvariety
in the Grassmanian, and
$T^*$ the subset of Grassmanian of lines in $\CP_4^*$ whose
elements are lines which are dual to 2-planes 
belonging to $T$.
This is a Zariski-open subset of the last Grassmanian.
This finishes preliminary considerations.

Since any irreducible component of the subvarieties $A$ and $\ol D_1,\dots, \ol D_k$ are at most 2-dimensional, by a dimensional reason,
there exists a line $l\subset\CP_4^*$ 
which does not intersect any of these subvarieties.
Moreover, the line $l$ can be taken from the subset $T^*$
since $T^*$ is Zariski-open in 
the Grassmanian of lines in $\CP_4^*$.
Let $P\subset\CP_4$ be the 2-plane which is dual to $l$, and $f':S'\to\CP_1
\simeq l$  the elliptic fibration
determined by the 2-plane $P$ as in \eqref{fprime}.
Since $l\in T^*$,  the intersection $S\cap P$ is transversal,
and in particular we have $p_i\not\in P$ for any $i=1,\dots, k$.
For $i=1,\dots, k$, we put $H_i:=l\cap p_i^*$.
In other words, $H_i$ is a hyperplane in $\CP_4$ spanned by $P$ and the point $p_i$.
Since $l\cap \ol D_i=\emptyset$, 
from \eqref{olD} we have $l\cap P_{ij}=\emptyset$ for any 
$j\neq i$.
This means  $p_j\not\in H_i$ for any $j\neq i$.
Hence, for any $i=1,\dots, k$, on the section $S\cap H_i$
there is no singular point of $S$ other than the point $p_i\in S$.
This means the property (i) in the lemma.
Moreover, the assumption $l\cap \ol D_i=\emptyset$
means that the section $S\cap H_i$ does not have
a singularity other than the point $p_i$.
This means that the fibration $f'$
satisfies the property (ii) in the lemma.

Next, let $H\in l$ be a hyperplane on which no singularity of $S$ belongs.
Because a hyperplane section $S_H$ is singular
only when $H$ belongs to the subvariety $\ms S$, 
by \eqref{Sdp}
the fiber $(f')\inv(H)\simeq S_H$ is singular
only when $H\in S^*$ or $H= H_i$
for some $i=1,\dots,k$.
From the choice of $H$, the latter cannot occur.
Moreover, since we have chosen a line $l$ which satisfies $l\cap A=\emptyset$, 
we have $H\not\in A=S^*\minus W_0$ and therefore
if $H\in S^*$ then $H\in W_0$.
Therefore, 
the section $S_H$ has a unique node as its all singularities and it is a smooth point of $S$.
In particular the fiber $(f')\inv(H)$ is of type I$_1$
and no singularity of $S'$ belongs to the same fiber.
This means that the fibration $f'$ satisfies the property (iii) in the lemma.
\proofend

\medskip
By taking the minimal resolutions 
of all singularities for the elliptic surface $S'$
which satisfies the three properties in 
the lemma,
we immediately obtain the following
\begin{corollary}\label{cor:elf}
Let $S\subset\CP_4$ be any Segre quartic surface
and $p_1,\dots, p_k$ the singularities of $S$.
Then there exists a 2-plane $P\subset\CP_4$ 
which intersects $S$ transversally at four points and 
for which the induced elliptic fibration
$\tilde f':\tilde S'\to\CP_1$ in \eqref{fprimet} satsisfies
the following properties.
\begin{itemize}
\item[(i)]
If $X_i$ is the dual graph of 
the exceptional curves of the singularity $p_i$,
then the dual graph of the singular fiber which includes
the exceptional curves of that singularity
is the extended Dynkin diagram $\tilde X_i$.
\item[(ii)]
All other singular fibers of $\tilde f'$ are
of type $I_1$.
\end{itemize}
\end{corollary}

Now we are able to prove the class formula for Segre surfaces.

\begin{theorem}\label{thm:cf}
Let $S\subset\CP_4$ 
and $X_1,\dots, X_k$ be as in the previous
corollary,
and $e_1,\dots, e_k$ the topological Euler characteristics
of the singular fibers of the elliptic fibration \eqref{fprimet},
so that the dual graphs of the fibers are of type $\tilde X_1,\dots,\tilde X_k$
respectively.
Then we have the formula
\begin{align}\label{cf1}
\deg S^* = 12 - (e_1+ \dots +e_k).
\end{align}
In particular we have $\deg S^*=12$ if 
the Segre surface $S$ is smooth.
\end{theorem}

\proof
We take a 2-plane $P\subset\CP_4$  as in the previous corollary
and let
$\tilde f':\tilde S'\to\CP_1$ be the associated elliptic fibration.
Then the set of critical values of $\tilde f'$
consists of hyperplanes $H_i = l\cap p_i^*$, $1\le i\le k$,
as well as the points corresponding to the singular fibers of type I$_1$.
The singular fibers over the former kind of critical values
are reducible since each of them includes the exceptional curves of the minimal resolution
of the singularity as components.
Hence the two kinds of critical values do not have a common point.
Moreover, we have $H_i\not\in S^*$ since  if not,
we have $H_i\in S^*$ but from our choice we have $l\cap A=\emptyset$ 
and therefore the hyperplane $H_i$ would belong to $W_0$,
which contradicts $p_i^*\cap W_0 = \emptyset$.
On the other hand, if $H\subset\CP_4$ is a hyperplane
for which $S\cap H$ corresponds to a singular fiber of type I$_1$,
the node of this singular fiber is a smooth point of $S$.
This means that $H\in W_0$.
These imply $l\cap S^*=l\cap W_0$, and that points belonging to the intersection $l\cap W_0$ 
are in one-to-one correspondence with the singular fibers of type I$_1$ of $\tilde f'$.

Since the elliptic surface $\tilde S'$ is rational and relatively minimal as seen before,
the topological Euler
characteristic of $\tilde S'$ is 12.
Hence, from the additivity of the topological Euler characteristic
to the elliptic fibration $\tilde f':\tilde S'\to\CP_1$ , we obtain that
the number of singular fibers which are of 
type I$_1$ is exactly $12 - (e_1+\dots + e_k)$.
From the conclusion in the first paragraph,
this directly implies   the desired equality \eqref{cf1}.
\proofend

\subsection{The Segre symbol}\label{ss:segs}
As is already mentioned, any Segre quartic surface
is a complete intersection of two quadrics in $\CP_4$. 
In general, complete intersections of two quadrics in $\CP_n$
for arbitrary $n$
can be systematically investigated 
by using so called the Segre symbol.
In this section,  following \cite[Chapter\,XIII,
Section\,10]{HP} and 
\cite[Section 8.6]{Dol},
We first recall Segre symbol
and next present a classification of Segre surfaces
in terms of the Segre symbol.

Let $X_0,\dots,X_n$ be homogeneous coordinates
on $\CP_n$.
Any quadric in $\CP_n$ is defined as the zero locus
of a quadratic form in $X_0,\dots,X_n$,
and  quadratic forms are in one-to-one correspondence with
symmetric matrices of size $(n+1)\times (n+1)$ in a standard way.
Let $Y$ be a complete intersection of two quadrics in $\CP_n$.
We write $\ms Q$ for the pencil of quadrics  generated by equations of $Y$.
Let $U$ and $V$ be symmetric matrices
which correspond to distinct two elements of $\ms Q$.
Assume that some member (and hence a generic member) of $\ms Q$ is smooth.
Then we may suppose $|V|\neq 0$.
If $\lmd$ is an indeterminate,
the determinant $|U-\lmd V|$ is a polynomial 
in $\lmd$ whose degree is precisely $n+1$.
Let
$$(\lmd - \aaa_1)^{e_1},\,\,
(\lmd - \aaa_2)^{e_2},\,\,\dots,\,\, (\lmd - \aaa_s)^{e_s}$$
be all elementary divisors of the matrix $U-\lmd V$,
so that $\aaa_1,\dots, \aaa_s$ are roots of 
the equation $|U-\lmd V|=0$.
Unlike the presentation in \cite{Dol} we do not assume $\aaa_i\neq \aaa_j$
for $i\neq j$, but we put indices for 
the roots $\aaa_i$ 
in such a way that the same roots are 
adjacent in the sense that 
 $\aaa_i=\aaa_j$ for some $i<j$ implies $\aaa_i = \aaa_{i+1} = \aaa_{i+2} = \dots = \aaa_j$.
From $|V|\neq 0$, 
we have
$$
0<s\le n+1
\qandq
\sum_{i=1}^s e_i = n+1.
$$
For a positive integer $e$,
we define two $e\times e$ matrices by 
\begin{align*}
P_e(\aaa) =
\begin{pmatrix}
0 & 0 &\dots & 0 & 0 & \aaa\\
0 & 0 & \dots & 0 &\aaa & 1\\
0 & 0 &\dots & \aaa & 1 & 0\\
\vdots &\vdots&\reflectbox{$\ddots$}&\vdots&\vdots&\vdots\\
0 & \aaa & \cdots & 0 & 0 & 0\\
\aaa & 1& \dots & 0 & 0 & 0
\end{pmatrix},
\quad
Q_e =
\begin{pmatrix}
0 & 0 &\dots & 0 & 0 & 1\\
0 & 0 & \dots & 0 &1 & 0\\
0 & 0 &\dots & 1 & 0 & 0\\
\vdots &\vdots&\reflectbox{$\ddots$}&\vdots&\vdots&\vdots\\
0 & 1 & \dots & 0 & 0 & 0\\
1 & 0& \dots & 0 & 0 & 0
\end{pmatrix},
\end{align*}
where $\aaa$ is any complex number.
Note that the elementary divisor of 
$P_e(\aaa) - \lmd Q_e$ is $(\lmd - \aaa)^{e}$.
The matrices $U$ and $V$ can be
simultaneously normalized in the sense of 
the theory of quadratic forms respectively
 into the matrices
$$
\begin{pmatrix}
 P_{e_1}(\aaa_1) & 0 & \dots &0\\
 0 & P_{e_2}(\aaa_2) & \dots & 0\\
 \vdots & \vdots & \ddots & \vdots\\
 0 & 0 & \dots & P_{e_s}(\aaa_s)
\end{pmatrix},
\quad
\begin{pmatrix}
 Q_{e_1} & 0 & \dots &0\\
 0 & Q_{e_2} & \dots & 0\\
 \vdots & \vdots & \ddots & \vdots\\
 0 & 0 & \dots & Q_{e_s}
\end{pmatrix}.
$$
Under this situation, the {\em Segre symbol}\, for the complete intersection $Y$ is given by
$$
[e_1e_2\dots e_s],
$$
with the exception that if 
some of the roots $\aaa_i$ are equal,
namely if 
$\aaa_i = \aaa_{i+1} = \aaa_{i+2} = \dots = \aaa_j$
for some $i<j$ and if all other roots are different from this common number, then the
sequence of  the entries $e_ie_{i+1}e_{i+2}\dots e_j$ is enclosed by round brackets.

For example, if $n=4$ and the Segre symbol is 
$[2111]$, the elementary divisors
of the matrix $U-\lmd V$ are $(\lmd - \aaa_1)^2,
\lmd - \aaa_2,\lmd - \aaa_3$
and $\lmd - \aaa_4$
for some 
distinct numbers $\aaa_1,\dots, \aaa_4$,
and we have $e_1 =2, e_2 = e_3 = e_4
=1$.
Hence the symmetric matrices $U$ and $V$ are 
simultaneously normalized respectively
 into the matrices
$$
\begin{pmatrix}
 0     & \aaa_1  & 0      & 0 & 0\\
\aaa_1 & 1       & 0      & 0 & 0\\
0      & 0       & \aaa_2 & 0 & 0\\
0      & 0       & 0      & \aaa_3 & 0\\ 
0      & 0       & 0      & 0 & \aaa_4
\end{pmatrix},
\quad
\begin{pmatrix}
 0     & 1       & 0      & 0 & 0\\
1      & 0       & 0      & 0 & 0\\
0      & 0       & 1      & 0 & 0\\
0      & 0       & 0      & 1 & 0\\ 
0      & 0       & 0      & 0 & 1
\end{pmatrix}.
$$
Therefore, the complete intersection 
$Y$ is defined by the equations
$$
2\aaa_1 X_0 X_1 +  X_1^2 + \aaa_2 X_2 ^ 2 
+\aaa_3 X_3 ^ 2 + \aaa_4 X_4 ^ 2=
2X_0X_1 + X_2^2 + X_3^2 + X_4^2 =0.
$$
As the second example, if $n=4$ and the Segre symbol is $[32]$,
all elementary divisors
of the matrix $U-\lmd V$ are $(\lmd - \aaa_1)^3$
and $(\lmd - \aaa_2)^2$ for some 
distinct numbers $\aaa_1$ and $\aaa_2$,
and we have $e_1 =3, e_2 = 2$.
The symmetric matrices $U$ and $V$ are 
simultaneously normalized respectively
 into the matrices
$$
\begin{pmatrix}
 0     & 0       & \aaa_1 & 0 & 0\\
 0     & \aaa_1  & 1      & 0 & 0\\
\aaa_1 & 1       & 0      & 0 & 0\\
0      & 0       & 0      & 0 & \aaa_2\\ 
0      & 0       & 0      & \aaa_2 & 1
\end{pmatrix},
\quad
\begin{pmatrix}
 0     & 0       & 1 & 0 & 0\\
 0     & 1  & 0      & 0 & 0\\
1 &     0   & 0      & 0 & 0\\
0      & 0       & 0      & 0 & 1\\ 
0      & 0       & 0      & 1 & 0
\end{pmatrix}.
$$
Therefore, the complete intersection 
$Y$ is defined by 
$$
2\aaa_1 X_0 X_2 + \aaa_1 X_1^2 + 2 X_1 X_2
+ 2\aaa_2 X_3 X_4 + X_4^2 = 
2X_0X_2 + X_1^2 + 2X_3 X_4 = 0.
$$
Instead, if the Segre symbol of $Y$ is $[(32)]$,
we just need to let $\aaa_1 = \aaa_2$ in 
this argument.
But the case $[(32)]$ does not give a Segre surface
since  it turns out from the normalized equations 
that it is a cone over a quartic curve.

A complete intersection of two quadrics in $\CP_4$
is smooth iff its Segre symbol is $[11111]$.
A list of Segre symbols for complete intersections
of two quadrics in $\CP_4$
which define Segre quartic surfaces
and whose pencil $\ms Q$ of quadrics contains
a smooth member
is given in Dolgachev's book \cite[p.\,398]{Dol},
and they consist of 16 kinds as in the following list in terms of the number of 
distinct roots of the equation $|U-\lmd V| = 0$,
or equivalently the number of singular members
of the pencil $\ms Q$.

\begin{itemize}
\item $[11111]$
\item $[2111]\quad [(11)111]$
\item $[(11)(11)1]\quad [(11)21]
\quad [311] \quad [221] \quad [(12)11]$
\item $[41]\quad [(31)1]\quad [3(11)]
\quad [32]\quad [(12)2]\quad [(12)(11)]$
\item $[5]\quad [(41)]$
\end{itemize}

The above assumption on smoothness for 
members of the pencil $\ms Q$ is in effect not necessary  (i.e.\,the above list covers all Segre quartic surfaces) by the following

\begin{proposition}\label{prop:cmp}
If the pencil $\ms Q$ of quadrics in $\CP_4$
does not have a smooth member,
the complete intersection defined by $\ms Q$ 
is not a Segre quartic surface.
\end{proposition}

\proof
We use the notations and results in \cite[Chap.\,XIII,
Sect.\,10 \& 11]{HP}.
Under the hypothesis on $\ms Q$,
we have $|V|=0$.
We again let
$(\lmd - \aaa_1)^{e_1},\dots, (\lmd - \aaa_s)^{e_s}$
be the elementary divisors of the matrix $U-\lmd V$,
where in the present situation $s=0$
if $|U-\lmd V|=0$.
There exist  integers $k>0$ and $r_0\ge 0$
 as well as 
a collection of integers $r_1,\dots, r_k>0$
which are determined from $U$ and $V$,
such that the relation
\begin{align}\label{rela}
\sum_{i=1}^s e_i = n+1 - 2\sum_{i=1}^k r_i - k - r_0
\end{align}
holds (\cite[p.290]{HP}).
The integer
$r_0$ is the number of variables among $X_0,\dots, X_n$ which do not appear in any equations of 
quadrics in $\ms Q$, and since Segre surfaces are not
a cone we have $r_0=0$.
Even so, since we have $n=4$ (as we are considering complete
intersections in $\CP_4$), if $k>1$ would hold, RHS of \eqref{rela} would be negative which is a contradiction.
So we have $k=1$, and the
relation \eqref{rela} becomes
\begin{align}\label{rela2}
\sum_{i=1}^s e_i = 4 - 2r_1.
\end{align}
Since LHS of this equation is non-negative, we have $r_1 = 1,2$.

If $r_1=1$, we have $\sum_{i=1}^s e_i=2$ and 
the possibilities for all elementary divisors of $U-\lmd V$ are
$$
\begin{matrix}
s=1 &(\therefore e_1 = 2)       &\Ra &(\lmd - \aaa_1)^2,\\
s=2 &(\therefore e_1 = e_2 = 1) &\Ra &(\lmd - \aaa_1),\,\,(\lmd - \aaa_2).
\end{matrix}
$$
Segre symbols of these cases are respectively
$[2;1]$ if $s=1$,
and $[11;1]$ or $[(11);1]$ if $s=2$ 
according as $\aaa_1\neq\aaa_2$ or $\aaa_1=\aaa_2$ respectively.
(The numbers put after the semicolon are $r_1,\dots, r_k$ in general.)
The normalized equations for these complete intersections are
(see \cite[p.294]{HP})
\begin{align*}
[2;1] &\Ra 2X_0X_1 + 2\aaa_1 X_3X_4 + X_4^2 = 
2X_1X_2 + 2X_3X_4 = 0,\\
[11;1] &\Ra 2X_0X_1 + \aaa_1 X_3^2 + \aaa_2 X_4^2=
2X_1X_2+ X_3^2 + X_4^2 = 0,\\
[(11);1] &\Ra 2X_0X_1 + \aaa_1 X_3^2 + \aaa_1 X_4^2=
2X_1X_2+ X_3^2 + X_4^2 = 0,
\end{align*}
where $\aaa_1\neq\aaa_2$ for the middle one.
The first one is reducible since it contains 
the 2-plane $X_1 = X_4=0$.
The last one is also reducible since it contains
the 2-plane $X_1 = X_3+iX_4=0$
(and also the plane $X_1 = X_3-iX_4=0$).
So these are not Segre surfaces.
For the middle one it can be readily seen that 
the surface has singularities along
the line $X_1 = X_3 = X_4=0$.
Since any singularity of a Segre surface is isolated, this means that
the surface is not a Segre quartic surface.

Finally, if $r_1 = 2$,
the Segre symbol of the complete intersection is 
$[;2]$, and the  normalized  equations are
$$
2X_0X_1 + 2X_2X_3 = 2X_1X_2 + 2X_3X_4=0.
$$
This is again reducible since it contains the plane 
$X_1 = X_3 = 0$, and hence 
does not give a Segre surface.
Therefore, $|V|\neq 0$.
\proofend

\subsection{Double covering structures on Segre quartic surfaces}
\label{ss:dcS}
As explained in the previous subsection,
the Segre symbol is useful for a systematic study
of complete intersections of two quadrics.
Especially, it provides us a pair of normalized forms  
for quadratic polynomials which define the surface.
In this subsection by using the normalized equations,
we see that
most Segre quartic surfaces have a structure
of a double covering over an irreducible
quadric in $\CP_3$.
We will use such a structure to determine
the types of singularities for most Segre surfaces.
We note that in \cite[p.401]{Dol} the types of all singularities of any Segre surfaces are presented without a proof.
The double covering structure will also be used 
in the next subsection 
to investigate complementary divisors
of the associated Einstein-Weyl spaces.

We begin with an easy observation
which is useful for our purpose.

\begin{proposition}\label{prop:SS}
Suppose that the symbol of a Segre quartic surface
$S$ has at least one `1' not living in a pair of round brackets.
Then by a generating projection from a point not lying on the surface,
$S$ has a structure of a finite double covering 
over a smooth quadric in $\CP_3$,
whose branch curve $B$ is a complete intersection
of the smooth quadric and another quadric.
Further, the Segre symbol of the curve $B$
is obtained from that of $S$ by just removing 
the `1'.
\end{proposition}

For example, if the Segre symbol of 
a Segre surface $S$ is $[1112]$,
since this contains (indeed three) `1' and
each of them do not live in a pair of round brackets,
$S$ is a double cover over a smooth quadric
surface whose branch is a complete
intersection with another quadric,
and the Segre symbol of the branch curve is $[112]$.
On the contrary, if 
the Segre symbol of 
a Segre surface $S$ is $[(11)(12)]$ for instance,
then since all `1' in the symbol are living in 
a pair of round brackets,
we cannot apply the proposition to such an $S$.


\medskip\noindent
{\em Proof of Proposition \ref{prop:SS}.}
The effect of changes of the order of the entry numbers in Segre symbol is just the exchanges for
the variables
in homogeneous coordinates on a projective space.
Hence, we may suppose that 
the `1' 
corresponds to the first variable $X_0$.
Since the `1' is not living in a pair of round brackets,
the corresponding elementary divisor is 
of the form $\lmd - \aaa_1$ for some $\aaa_1\in\CC$.
Let $F$ and $G$ be a pair of normalized 
quadratic polynomials of $S$ in these coordinates.
From the above elementary divisor,
monomials in $F$ and $G$ which contain the variable
$X_0$ are (a constant multiple of) $X_0^2$ only,
and we have
\begin{align}\label{qeq1}
F = \aaa_1 X_0^2 + f,\,\,\,
G = X_0^2 + g,
\end{align}
where the residual quadratic polynomials $f$ and $g$ do not contain $X_0$.
Let $\pi:\CP_4\to\CP_3$ be the generic projection
from the point $(1,0,0,0,0)$.
Namely $\pi$ is the map which drops the coordinate
$X_0$.
Eliminating $X_0^2$ from the equations $F=G=0$, we obtain
\begin{align}\label{qeq2}
f - \aaa_1 g = 0.
\end{align}
This is a quadratic equation in $X_1,X_2,X_3,X_4$,
and the corresponding symmetric $4\times 4$ matrix
is readily seen to be of full rank (i.e.\,rank 4)
by using that $\aaa_1$ is different from other roots
appearing in the polynomial $f$.
Hence if $Q\subset\CP_3$ is the quadric which is 
defined by \eqref{qeq2},
it is smooth.
From \eqref{qeq1}, the generating point of the projection 
$\pi$ is not on $S$,
and the restriction $\pi|_S$ provides
$S$ a structure of finite double cover over $Q$,
with the branch divisor $B$ being defined by 
$f-\aaa_1 g = g=0$.
Hence the branch divisor is a complete intersection of two quadrics in $\CP_3$.
The last equations are equivalent to $f=g=0$,
so we have $B = \{f=g=0\}$.
Moreover, obviously the two polynomials $f$ and $g$ are already in  normalized forms,
and the Segre symbol of the complete intersection
$B$ is exactly the one given in
the proposition.
\proofend

\begin{figure}
\includegraphics{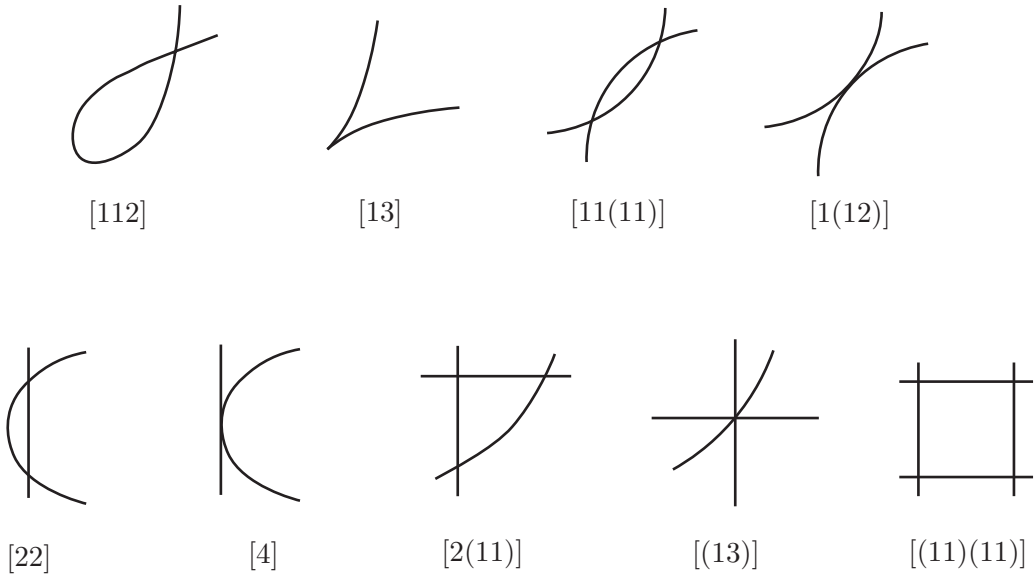}
\caption{
The branch curves $B$
}
\label{fig:quartics}
\end{figure}

\medskip
When the symbol of a Segre surface
is as in the proposition,
principal structure of the surface may be  read off from the structure of 
the branch curve $B$ in the smooth quadric $Q$.
In turn structure of the curve $B$ is known
from the Segre symbol of $B$,
because the symbol provides normalized forms
of the quadratic equations of $B$,
and it is not difficult to obtain concrete
structure from the equations,
by projecting $B$ to $\CP_2$ from a point.
The cases $[1111],[112]$ and $[13]$ will be
treated in Propositions \ref{prop:e1}--\ref{prop:e3}
for other purpose.
Below we present a list of structure of $B$ which actually arise from Segre surfaces as 
in Proposition \ref{prop:SS}.
These can also be found in \cite[pp.\,305--8]{HP}
where arbitrary complete intersections of two quadrics in $\CP_3$
are treated. See also Figure \ref{fig:quartics} for these curves.

\medskip\noindent
$[1111]\Ra$ smooth elliptic curve,
\,
$[112]\Ra$ 1-nodal rational curve,
\,
$[13]\Ra$ 1-cuspidal rational curve, 
\,
$[11(11)]\Ra$ two conics intersecting transversally at two points,
\,
$[1(12)]\Ra$ two conics touching at one point,
\, 
$[22]\Ra$ one line and one rational normal curve
intersecting transversally at two points,
\,
$[4]\Ra$ one line and one rational normal curve
touching at one point,
\,
$[2(11)]\Ra$ two lines and one conic, forming a `triangle',
\,
$[(13)]\Ra$ two lines and one conic, sharing one point,
\,
$[(11)(11)]\Ra$ a `square' of four lines.

From these, the types of all singularities of Segre 
surfaces whose symbols satisfy the property
in Proposition \ref{prop:SS} are obtained.
We present them in the second column of Table \ref{table1}.
Then by using Theorem \ref{thm:cf}, 
the classes of these Segre surfaces are obtained, 
and they are listed in the third column in Table \ref{table1}.
The forth column of the table presents the numbers of `1' in the symbols
which are not living in a pair of round brackets.
These are  in one-to-one correspondence with
generating projections from a point which induces the double
covering map to the smooth quadric surface as in Proposition \ref{prop:SS}.
We will show in Section \ref{ss:sil} that the dual quadric $Q^*$
of the smooth quadric $Q$ is always contained in the dual variety $S^*$,
and this is why we are writing $Q^*$ in the table.
The numbers of lines on $S$ are listed in the fifth column,
and they are taken from \cite[p.\,401]{Dol}.
Among them the number of lines on $S$
which do not pass any singularity of $S$ are listed in 
the sixth column.
In Section \ref{ss:sil} we will also show that
the 2-plane in $\CP_4^*$ which is formed by 
hyperplanes in $\CP_4$ that contain such a line
is always contained in $S^*$.

\begin{table}[h]
\begin{tabular}{|c|c|c|c|c|c|c|c|}
\hline
Segre symbol & 
$\Sing S$ & $\deg S^*$  & $\# \{Q^*\!\subset \!S^*\}$ &
$\#\{$lines\,$\subset \!S\}$ & $\#\{\CP_2\!\subset \!S^*\}$ &
$\Aut_e S$\\
\hline\hline
[11111] & none & 12 & 5 & 16  & 16 
 & \{\id\}\\
\hline\hline
[1112] & ${{\rm A}}_1$ & 8 & 3 & 12 & 8 
 & \{\id\} \\
\hline
[111(11)] & $2{{\rm A}}_1$ & 8 & 3 &  $8$ & 0 
&$\CC^*$ \\
\hline\hline
[12(11)] & $3{{\rm A}}_1$ & 6 & 1 & 6 & 0  & 
 $\CC^*$\\
\hline
[1(11)(11)] & $4{{\rm A}}_1$ & 4 & 1 & 4 & 0 
& $\CC^*\times\CC^*$ \\
\hline
$[113]$ & ${{\rm A}}_2$ & 9 & 2 & 8 & 4  
& $\{\id\}$
\\
\hline
$[122]$ & $2{{\rm A}}_1$ & 8  & 1 & 9 & 4
& $\{\id\}$\\
\hline
$[11(12)]$ & ${{\rm A}}_3$ & 8  & 2  & 4 & 0 
& $\{\id\}$ \\
\hline\hline
$[14]$ & ${{\rm A}}_3$  & 8 & 1 &  5 & 2 
 & $\CC^*$\\
\hline
$[1(13)]$ & ${{\rm D}}_4$ & 6  & 1  & 2 & 0 
 & $\CC^*$\,or $\{\id\}$\\
\hline
\end{tabular}
\bigskip
\caption{Segre quartic surfaces which are realizable as a double cover over a smooth quadric surface}
\label{table1}
\end{table}

Next, we discuss the case where the Segre symbol
contains `1' but where it lives in a pair of round brackets.
Then the situation is little different from 
that in the previous proposition,
but it happens to be still simple.

\begin{proposition}\label{prop:SC}
Suppose that the Segre symbol of a Segre quartic surface $S$ includes `1' 
which is enclosed by a pair of round brackets.
Then by a generating projection from a point not lying on the surface,
$S$ has a structure of a finite double covering 
over the cone over an irreduble conic,
and the branch curve is a complete intersection
of the cone with another quadric.
Further, the Segre symbol of the branch curve
is obtained from that of $S$ by 
removing the `1' 
and next removing the pair of round brackets
in which the `1' lives.
\end{proposition}

For example, the symbol
$[(11)12]$ (resp.\,$[(14)]$) enjoys the assumption
of the proposition, and the symbol of the branch curve
on the cone is $[112]$ (resp.\,$[4]$).
Combining with Proposition \ref{prop:SS},
this means that if the symbol has two `1',
and if one of them is included in 
a pair of round brackets while the other one
is not
(like $[(11)12]$),
then $S$ is a double cover of a smooth quadric
and also is a double cover of the cone over an irreducible conic.

\medskip\noindent
{\em Proof of Proposition \ref{prop:SC}.}
From the list for the Segre symbols of all Segre surfaces given right before Proposition \ref{prop:cmp},
the Segre symbols in the situation of the present proposition
are of the forms $[(11)\bm e], [(12)\bm e']$
for some $\bm e$ and $\bm e'$,
or exactly one of $[(13)1]$ and $[(14)]$.
%
%
We verify the proposition on case-by-case basis
depending on these distinctions.

For the first case (in which the symbol is of the form $[(11)\bm e]$),
$\bm e$ is one of 
\begin{align}\label{Sgr3}
111, \,\,12, \,\,3, \,\,(11)1, \,\,(12).
\end{align}
The pair of normalized symmetric matrices which define
the surface $S$ are of the forms
\begin{equation}\def\arraystretch{1}
U=    \left(\begin{blockarray}{cc|cc}
        \aaa_1 & 0 & \BAmulticolumn{2}{c}{\multirow{2}{*}{\Large$0$}}
        \\
        0 & \aaa_1 &  
        \\
        \cline{1-4}
        \BAmulticolumn{2}{c|}{\multirow{2}{*}{\Large$0$}}&
        \BAmulticolumn{2}{c}{\multirow{2}{*}{$U'$}}
    \end{blockarray}
    \right),
    \quad
    V=    \left(\begin{blockarray}{cc|cc}
        1 & 0 & \BAmulticolumn{2}{c}{\multirow{2}{*}{\Large$0$}}
        \\
        0 & 1 &  
        \\
        \cline{1-4}
        \BAmulticolumn{2}{c|}{\multirow{2}{*}{\Large$0$}}&
        \BAmulticolumn{2}{c}{\multirow{2}{*}{$V'$}}
    \end{blockarray}
    \right),
\end{equation}
where $U'$ and $V'$ are $2\times 2$ matrices 
which are determined from $\bm e$.
In terms of these symmetric matrices,
an elimination of  the variable $X_0$
(or equivalently $X_1$) from the pair of the quadratic equations
corresponds to considering the matrix
$U - \aaa_1 V$, 
and we obtain a quadratic polynomial in $X_2,X_3,X_4$.
The  $3\times 3$ symmetric matrix corresponding to this polynomial
is exactly $U' - \aaa_1 V'$, and
it is not difficult to see from the list \eqref{Sgr3}
using $\aaa_i\neq \aaa_j$ for $i\neq j$,
that this matrix is of full rank.
So the quadratic polynomial defines an irreducible conic in $\CP_2$.
The point $(1,0,0,0,0)$ is again not on $S$,
and if $\pi:\CP_4\to \CP_3$ is the projection from this point,
the image $\pi(S)\subset\CP_3$ is contained in the cone over
the last conic.
The restriction $\pi|_S:S\to\pi(S)$ is of degree-two over
the image, so $\pi(S)$ is 2-dimensional.
Hence we obtain that the image $\pi(S)$ is 
the cone over the conic.
Moreover, we find that
the branch divisor of $\pi|_S$ is 
a complete intersection of two quadrics in $\CP_3$
whose symmetric matrices are
\begin{align}\label{matp}
\begin{pmatrix}
0 & 0 \\
0 & U'-\aaa_1V'
\end{pmatrix}
\qandq
\begin{pmatrix}
1 & 0 \\
0 & V'
\end{pmatrix}.
\end{align}
Again, by using $\aaa_1\neq\aaa_i$ if $i>1$,
it is readily seen (again by case-by case checking
relying on \eqref{Sgr3}) that
this pair is already in normal forms,
and that the corresponding Segre symbol 
is exactly $[1\bm e]$, the first `1' being
corresponding to the $1\times 1$-components of the upper-left 
in the matrices \eqref{matp}.
Thus, the Segre symbol of the branch curve
of $\pi|_S$ is certainly as in the proposition.

Next, if the Segre symbol of $S$ is of the form 
$[(12)\bm e']$ for some  $\bm e'$,
all the possibilities for $\bm e'$ are $11, (11)$ and $2$.
Then we have
\begin{equation}\def\arraystretch{1}
U=    \left(\begin{blockarray}{ccc|cc}
        \aaa_1 & 0 & 0 & \BAmulticolumn{3}{c}{\multirow{3}{*}{\large$0$}} &
        \\
        0 & 0 & \aaa_1 &  &
        \\
        0 & \aaa_1 & 1 & &
        \\
        \cline{1-6}
        &\BAmulticolumn{2}{c|}{\multirow{3}{*}{\large$0$}} & &
        \BAmulticolumn{2}{c}{\multirow{2}{*}{$U'$}} 
    \end{blockarray}
    \right),
    \quad
V=    \left(\begin{blockarray}{ccc|cc}
        1 & 0 & 0 & \BAmulticolumn{3}{c}{\multirow{3}{*}{\large$0$}} &
        \\
        0 & 0 & 1 &  &
        \\
        0 & 1 & 0 & &
        \\
        \cline{1-6}
        &\BAmulticolumn{2}{c|}{\multirow{3}{*}{\large$0$}} & &
        \BAmulticolumn{2}{c}{\multirow{2}{*}{$V'$}} 
    \end{blockarray}
    \right),
\end{equation}
where $U'$ and $V'$ are two symmetric matrices that are determined
from $\bm e'$.
By the same procedure as in the last case of $[(11)\bm e]$,
we obtain that the projection $\pi:\CP_4\to\CP_3$ from the point
$(1,0,0,0,0)\not\in S$ induces a double covering map to the
cone over an irreducible conic, and the branch curve
is a complete intersection of two quadrics in $\CP_3$ defined by two symmetric matrices
\begin{align}\label{matq}
\left(\begin{blockarray}{cc|cc}
        0 & 0 & \BAmulticolumn{2}{c}{\multirow{2}{*}{\large$0$}}
        \\
        0 & 1 &  
        \\
        \cline{1-4}
        \BAmulticolumn{2}{c|}{\multirow{2}{*}{\large$0$}}&
        \BAmulticolumn{2}{c}{\multirow{2}{*}{$U' -\aaa_1 V'$}}
    \end{blockarray}
    \right)
\qandq
\left(\begin{blockarray}{cc|cc}
        0 & 1 & \BAmulticolumn{2}{c}{\multirow{2}{*}{\large$0$}}
        \\
        1 & 0 &  
        \\
        \cline{1-4}
        \BAmulticolumn{2}{c|}{\multirow{2}{*}{\large$0$}}&
        \BAmulticolumn{2}{c}{\multirow{2}{*}{$V'$}}
    \end{blockarray}
    \right).
\end{align}
Again, by using $\aaa_1\neq \aaa_i$ for any $i>1$,
we can see that this pair is already in normal forms,
and the corresponding Segre symbol is $[2\bm e']$,
where the first `2' being
corresponding to the $2\times 2$ matrices in the upper-left 
in the matrices \eqref{matq}.
Thus,  the Segre symbol of the branch curve
of $\pi|_S$ is again  as in the proposition.

The cases where
the Segre symbol of $S$ is
$[(13)1]$ or $[(14)]$
can be shown in the same manner,
and we omit them.
\proofend

\medskip

When Proposition \ref{prop:SC} can be applied to 
a Segre surface $S$, since the double covering map is finite,
$S$ always has singularity over
the vertex of the cone.
Let $v$ be the vertex of the cone
and $B$ the branch divisor on the cone.
If $v\not\in B$, then
the surface $S$ has two ${{\rm A}}_1$-singularities over $v$.
If $v\in B$, then 
the surface $S$ has exactly one singularity over $v$.
Since the curve $B$ can have a singularity at a smooth point 
of the cone, $S$ can have other singularity in general,
but again the types of them are known from that of the branch curve $B$.
When $v\in B$, the type of the singularity of $S$ over $v$ 
is known by transforming $B$
through blowing up at $v$, and noticing that
the exceptional curve of the blowup is always included in 
the branch locus of the new double covering.
This way the types of all singularities
of $S$ are again known, and 
we display them in Table \ref{table2}. 
Again, the classes of the surfaces are obtained from Theorem \ref{thm:cf}, 
and we list them in the third column.
Positioning of the vertex $v$ and the branch curve $B$ is shown in the fifth column.
We note that it is possible to show that any line on $S$
passes a singularity of $S$ for any Segre surface $S$
which appears in Table \ref{table2},
and therefore we obtain no plane in the dual variety
$S^*$ in the situation
of Proposition \ref{prop:SC}.
Also, we note that $[(11)(12)]$ and $[(12)(11)]$ 
represent  Segre surfaces of the same kind,
and we are taking  different projections to the cone.

\begin{table}[h]
\begin{tabular}{|c|c|c|c|c|c|c|c|}
\hline
Segre symbol &  $\Sing S$ &
$\deg S^*$  & $\#\{$lines $\subset S\}$ & the vertex $v$
&$\Aut_eS$\\
\hline\hline
$[(11)111]^*$ & $2{{\rm A}}_1$ & 8 
& 8 & $v\not\in B$ &$\CC^*$ \\
\hline\hline
$[(12)11]^*$ & ${{\rm A}}_3$ & 8 & 4 &  $v$ is a node of $B$
& $\{\id\}$  \\
\hline
$[(11)12]^*$ & $3{{\rm A}}_1$ & 6  & 6 & $v\not\in B$ & $\CC^*$ \\
\hline
$[(11)(11)1]^*$ & $4{{\rm A}}_1$ & 4 & 4 & $v\not\in B$ 
& $\CC^*\times\CC^*$ \\
\hline\hline
$[(11)3]$ & $2{{\rm A}}_1+{{\rm A}}_2$ & 5 & 4 
& $v\not\in B$
 & $\CC^*\times\CC^*$ \\
\hline
$[(13)1]^*$ & ${{\rm D}}_4$ & 6  & 2  & $v$ is a cusp of $B$
  & $\CC^*$ or $\{\id\}$\\
\hline
$[(12)2]$ & ${{\rm A}}_1+{{\rm A}}_3$  & 6  & 3
&  $v$ is a node of $B$ & $\CC^*$\\
\hline
$[(11)(12)]$  & $2{{\rm A}}_1+{{\rm A}}_3$ & 4  & 2 &
 $v$ is a node of $B$ & $\CC^*\times\CC^*$ \\
\hline
$[(12)(11)]$& $2{{\rm A}}_1+{{\rm A}}_3$ & 4  & 2 & 
$v\not\in B$ & $\CC^*\times\CC^*$  \\
\hline\hline
$[(14)]$ & ${{\rm D}}_5$ & 5   & 1
& $\Sing B = \{v\}$ & $\CC^*$ or $\{\id\}$\\
\hline
\end{tabular}
\medskip
\caption{Segre quartic surfaces which are realizable as a double cover of the cone over an irreducible conic.
Symbols which appeared in Table \ref{table1} are associated by `$*$'. }
\label{table2}
\end{table}
Evidently at least one of Propositions \ref{prop:SS}
and \ref{prop:SC} is applied if the Segre symbol contains at least one `1'. 
Among 16 symbols for Segre quartic surfaces,
there are exactly two symbols which do not
contain any `1', and they are 
\begin{align}\label{}
[23] \qandq [5].
\end{align}
According to \cite[Table 8.6]{Dol}, the former surface
has exactly two singularities
and they are of types ${{\rm A}}_1$ and  ${{\rm A}}_2$,
and the latter surface has exactly one
singularity, which is of type ${{\rm A}}_4$.
For these Segre surfaces,
we cannot obtain a double covering structure
over a quadric by 
a generic projection from a point which does not
belong to $S$.
Similarly to the previous two tables, 
we display principal structures of these surfaces in Table \ref{table3}.

\begin{table}[h]
\begin{tabular}{|c|c|c|c|c|c|c|}
\hline
Segre symbol &  $\Sing S$ & $\deg S^*$ & \#$\{$lines\,$\subset S\}$ & 
\#$\{\CP_2\subset S^*\}$   & $\Aut_e S$ \\
\hline\hline
[23] & ${{\rm A}}_1 + {{\rm A}}_2$ & 7 & 6 & 3 &  $\CC^*$ \\
\hline\hline
[5] & ${{\rm A}}_4$ & 7& 3 &  1 &  $\CC^*$\\
\hline
\end{tabular}
\bigskip
\caption{Segre quartic surfaces which cannot
be realized as a double cover over a quadric surface}
\label{table3}
\end{table}


\subsection{Self-intersection loci of the dual varieties.}
\label{ss:sil}
Recall that if $S\subset\CP_4$ is a Segre quartic surface, the dual variety $S^*\subset\CP_4^*$ 
contains a Zariski-open subset $W_0$ parameterizing
1-nodal rational curves that are contained
in $S_{\reg}$, the smooth locus of $S$,
and $W_0$ is a complex 3-manifold which has
an Einstein-Weyl structure (see Definition \ref{def:Severi} and Proposition \ref{HNMT}).
The subset $W_0$ cannot coincide with the full set $S^*$ since 1-nodal rational curves can always be 
deformed in $S$ to a 2-nodal (and hence reducible) curve or a 1-cuspidal rational curve.
So the complement $S^*\minus W_0$ is always a non-empty subvariety in $S^*$.
This subvariety may be regarded as
the natural boundary set or a kind of conformal infinity of the complex 
Einstein-Weyl space $W_0$.
As in the introduction, we call 2-dimensional components 
of $S^*\minus W_0$ {\em divisors at infinity}.
In this subsection we find several divisors at infinity
in concrete forms,
and show that the dual varieties have self-intersection
along these divisors.

In order to find and describe divisors at infinity,
we first make use of lines lying on a Segre surface.
Let $S\subset\CP_4$ be a Segre quartic surface and  $l$  a line on $S$.
The number of lines on Segre surfaces are presented in 
Tables \ref{table1}--\ref{table3}.
We write $l^*$ for the set of hyperplanes in $\CP_4$ which 
contain the line $l$. So $l^*$ is a 2-plane in the dual space
$\CP_4^*$.
If $H\in l^*$, the line $l$ is contained in the
hyperplane section $S\cap H$, and
since $S$ is of degree four, the section is reducible.
This means that the section $S\cap H$ is not a minitwistor line.
So it does not belong to the open subset $W_0$ of $S^*$.
But we have 

\begin{proposition}\label{prop:line}
Suppose that the line $l\subset S$ does not pass any singularity
of $S$. Then the 2-plane $l^*\subset\CP_4^*$ is contained in 
the dual variety $S^*$. In particular, 
it is a divisor at infinity.
\end{proposition}

For the proof of the proposition, we first show

\begin{lemma}\label{lemma:lorc}
Let $S\subset\CP_4$ be a Segre quartic surface.
Then any 2-plane in $\CP_4$ does not 
contain a curve on $S$ whose degree in $\CP_4$ is greater than two.
\end{lemma}

\proof
Since any Segre surface is of degree four,
any curve on $S$ which is contained in a 2-plane is at most of degree four.
If there is a 2-plane $P\subset\CP_4$ 
such that $P\cap S$ contains a quartic curve,
then the pencil of hyperplanes in $\CP_4$ which contain the plane $P$
satisfies the property that for any hyperplane $H$ belonging to the
pencil the intersection $H\cap S$ is exactly the quartic curve.
Hence, $H\cap S$ is independent of the choice of $H$ in the pencil.
This contradicts that the restriction homomorphism
\begin{align}\label{isom66}
H^0 \big(\ms O_{\CP_4}(1)\big) \to
H^0\big(K_S\inv\big) 
\end{align}
is isomorphic.
So there is no 2-plane $P$ such that $P\cap S$ contains a quartic curve.
Next if there is a 2-plane $P\subset\CP_4$ 
such that the intersection $P\cap S$ contains a cubic curve,
the pencil of hyperplanes in $\CP_4$ which contain the plane $P$
satisfies the property that for any hyperplane $H$ belonging to the
pencil the intersection $H\cap S$ contains a line as the residual component
of the cubic curve in $P$.
Again from the natural isomorphism \eqref{isom66},
the residual line has to move as the hyperplane moves in the pencil.
This means that the surface $S$ contains infinite number of lines,
which is again a contradiction.
So there is no 2-plane $P$ such that $P\cap S$ contains a cubic curve.
Clearly these mean the assertion in the lemma.
\proofend

\medskip
The situation for conics on Segre surfaces is very different:

\begin{lemma}\label{lemma:cnc}
If a Segre surface $S\subset\CP_4$ contains an irreducible conic $C$
which does not pass any singularity of $S$,
we have $C^2 = 0$ on $S$.
\end{lemma}

\proof
Let $C_1$ be an irreducible conic on $S$
and suppose that it does not pass any singularity of $S$.
Let $P_1$ be the 2-plane which contains $C_1$.
Take a generic hyperplane $H\subset\CP_4$
such that $H$ contains the 2-plane $P_1$.
Then since $S$ is of degree four, the hyperplane section
$S\cap H$ is a quartic curve, and it contains the conic $C_1$.
Hence $S\cap H$ contains an irreducible conic
or two lines as the residual components.
But since there is a finite number of lines on $S$,
by choosing a generic $H$ containing $P$, we may suppose that 
the residual component of $S\cap H$ is an irreducible conic,
say $C_2$.
This determines a 2-plane $P_2$ which satisfies $C_2\subset P_2$.
Since $H|_S = C_1 + C_2$, the divisor $C_1+C_2$ is an anti-canonical curve on $S$.
Hence, we have $(C_1+C_2)^2 = 4$.
So we have
\begin{align}\label{2C1C2}
C_1^2 + C_2 ^2 + 2C_1C_2 = 4.
\end{align}
A priori this is an equation in $\QQ$, but 
since both $C_1$ and $C_1+C_2$ are Cartier divisors,
so is the curve $C_2$.
Hence not only $C_1^2$ but also $C_2^2$ and $C_1C_2$ are integers.
So \eqref{2C1C2} is an equation in $\ZZ$.
Since $C_1+C_2$ is an ample divisor, by Lefschetz  theorem
on hyperplane sections,
it is connected. Hence $C_1C_2>0$.
On the other hand, by adjunction formula and positivity 
of the bundle $K_S\inv$, we readily obtain $C_1^2\ge -1$
(and also $C_2^2\ge -1$).
Moreover, since the conic $C_2$ moves on $S$ by 
moving the hyperplane $H$ containing the plane $P_1$, 
we have $C_2^2\ge 0$.
Hence we have $C_1^2 + C_2^2 \ge -1$.
But by \eqref{2C1C2}, the number $C_1^2 + C_2^2$ is even,
so  $C_1^2 + C_2^2 \ge 0$.
These mean that  $C_1C_2$ is either $1$ or $2$.

If $C_1C_2=1$, from \eqref{2C1C2}, all possibilities for
the set $\{C_1^2, C_2^2\}$ of self-intersection numbers
are $\{2,0\}, \{1,1\}$ and $\{-1,3\}$.
We can readily see by standard calculations  using exact sequences
and rationality of $S$,
that in all these possibilities,
we have $\dim |C_1+C_2| = 5$.
This contradicts $\dim |K_S\inv|=4$.
Hence $C_1C_2 = 2$ has to hold.
From \eqref{2C1C2} this means $C_1^2 = C_2^2=0$.
In particular we have $C_1^2 =0$ on $S$.
This is the assertion of the lemma.
\proofend

\medskip
Another property on conics on Segre surfaces that will be used to prove Proposition \ref{prop:line} is the following

\begin{lemma}\label{lemma:cnc2}
Let $S\subset\CP_4$  be a Segre surface.
For each singularity of $S$,
conics on $S$ which pass that singularity constitute at most 1-dimensional
family.
\end{lemma}
\proof
Let $p_i$ be a singular point of $S$
and  $C_1\subset S$  an irreducible conic
which passes $p_i$.
Then again
for a generic hyperplane $H\subset\CP_4$ which 
contains the 2-plane that includes $C_1$,
the intersection $S\cap H$ is of the form
$C_1 + C_2$, where $C_2$ is another irreducible conic.
Let $S' \to S$ be the minimal resolution of the singularity $p_i$, and 
$C_1'$ and $C_2'$ the strict transforms of $C_1$
and $C_2$ into $S'$ respectively.
To show the lemma, it is enough to show $\dim|C'_1| \le 1$,
and since $C'_1\simeq C_1$ is a smooth rational curve,
this inequality is equivalent to $(C'_1)^2\le 0$.
So suppose $(C'_1)^2 > 0$.
Then the original conic $C_1$ moves on $S$ in at least a 2-dimensional family 
of conics that passes the singularity $p_i$.
To each conic in this family we obtain a pencil of residual conics on $S$ as above,
Therefore, we obtain a $(2+1=)$ 3-dimensional family of reducible anti-canonical curves
on $S$.
By a dimensional reason the parameter space of this family has to be
a component of the dual variety $S^*$, but this contradicts
the irreducibility of $S^*$.
Hence, we have  $(C'_1)^2\le 0$, which means the lemma as above.
\proofend

\medskip
By using these three lemmas, we next show

\begin{lemma}\label{lemma:scnt}
A generic secant $\ol{pq}$, $p,q\in S$,
of any Segre quartic surface $S\subset\CP_4$ has the properties that
it does not intersect any line on $S$, and that 
any 2-plane which contains $\ol{pq}$
does not include any curve on $S$.
\end{lemma}

\proof
By Lemmas \ref{lemma:cnc} and \ref{lemma:cnc2},
any component of the variety of conics in $\CP_4$ which are contained in $S$
is at most 1-dimensional.
A generic element $C$ of this variety 
defines a 2-plane $P\subset\CP_4$ by the condition $C\subset P$.
Hence the variety of conics on $S$ is identified with 
the variety of 2-planes in $\CP_4$ which contains 
a conic on $S$.
We then define $V$ to be the variety consisting of lines in $\CP_4$
which are contained in some 2-plane belonging to the last variety of 2-planes.
In other words, $V$ is the set of lines lying on the same 2-plane in $\CP_4$
as some conic on $S$.
Then obviously each component of $V$ is at most 3-dimensional.
Moreover, if $W$ denotes the variety formed by secants of the surface $S$,
then we have $\dim W = 4$, and moreover $V\subset W$.

Next let $K$ be the set of secants of $S$ which intersect some line in $S$.
Evidently $K$ is a 3-dimensional subvariety  of $W$.
By a dimensional reason, we may choose a secant $\ol{pq}\in W$
which does not belong to the subvariety $V\cup K$.
Then $\ol{pq}\cap l=\emptyset$ for any line $l\subset S$.
Let $P\subset\CP_4$ be any 2-plane containing $\ol{pq}$.
By Lemma \ref{lemma:lorc},  $P$ does not contain any curve on $S$
whose degree is greater than two.
Moreover from $\ol{pq}\not\in V$, the 2-plane $P$ does not contain
any conic on $S$.
\proofend

\medskip\noindent
{\em Proof of Proposition \ref{prop:line}.}
Let $l$ be a line on a Segre surface $S\subset\CP_4$.
Let $H$ be a generic hyperplane which contains $l$.
We first show that there exists a secant
of $S$ which is included in $H$ and which satisfies
the property in Lemma \ref{lemma:scnt}.
Since $S$ is quartic we have $S\cap H = l \cup C$,
where $C$ is a cubic curve.
By Lemma \ref{lemma:lorc}, $C$ is 
non-degenerate in the sense that it is not contained in a 2-plane.
Any secant of $S$ included in $H$
is a secant of the hyperplane section $S\cap H=l\cup C$.
Any secant of the cubic curve $C$ is of course such a secant,
and since $C$ is a curve they constitute 2-dimensional family of secants of $S$.
Because $C$ is non-degenerate as above,
generic two distinct secants in the family do not intersect.
If a generic secant of $S$ in $H$
does not satisfy the property in Lemma \ref{lemma:scnt},
then to each secant of $C$ we have a conic in $S$
which lies on the same 2-plane as the secant.
These mean that there exists a 2-dimensional family of conics in $S$,
which contradicts Lemmas \ref{lemma:cnc} and \ref{lemma:cnc2}.
Therefore, there is a secant $\ol{pq}$ of $S$
which is contained in $H$ and which satisfies the property 
in Lemma \ref{lemma:scnt}.
Obviously, we can choose such a secant which does not pass the line $l$.
Then the hyperplane $H$ is spanned by $l$ and $\ol{pq}$.

We fix any one of such secants $\ol{pq}\subset H$, and 
let $\varpi:\CP_4\to\CP_2$ be the projection from the line $\ol{pq}$.
Fibers of $\varpi$ are 2-planes which contain $\ol{pq}$.
We use the same letter $\varpi$ for the restriction of the projection
to the surface $S$.
Then $\varpi:S\to\CP_2$ has the two points $p$ and $q$
as the set of indeterminacy,
and if $S'$ is the blowing-up of $S$ at these two points,
then the composition $S'\to S\stackrel{\varpi}\to\CP_2$ is a morphism.
We write $\varpi'$ for this morphism.
This is the anti-canonical map from $S'$.
Since $S$ is of degree four, the morphism $\varpi'$ is 
of degree two.
Moreover, since there exists no 2-plane in $\CP_4$ which
contains a curve on $S$,
the degree-two morphism $\varpi':S'\to\CP_2$ does not 
contract any curve on $S'$.
Namely $\varpi'$ is a finite morphism.
Since $K_S^2 = \deg S = 4$, we have $K_{S'}^2 = 2$.
From these we readily obtain that the branch divisor
of $\varpi':S'\to\CP_2$ is a quartic curve.

Now assume that the line $l$ on $S$ does not to pass any singularity of $S$.
Then $l$ is a $(-1)$-curve on $S$.
Let $l'$ be the strict transform of $l$ into $S'$.
Then since the secant $\ol{pq}$ does not intersect $l$
from our choice,
$l'$ is still a $(-1)$-curve.
Hence we have $K_{S'}\inv.\,l' = 1$.
So the image $\varpi'(l')= \varpi(l)\subset\CP_2$ is a line,
and $l'$ is isomorphic to this line by $\varpi'$.
Further the pullback $(\varpi')\inv(\varpi'(l'))$ is an anti-canonical curve on
$S'$.
On the other hand, since the hyperplane $H$ contains the center
of the projection $\varpi$, the image $\varpi(H)$ is also a line,
and as $l\subset H$, we have $\varpi'(l') = \varpi(H)$.
If the line $\varpi'(l')$ would be contained in the branch quartic of $\varpi'$,
the branch curve is singular at some point on the line $\varpi'(l')$.
Since $\varpi'$ does not contract any curve on $S'$,
this means that there would exist a singular point of $S'$ on $l'$.
Hence $S$ would have a singularity at some point of the line $l$,
which contradicts the assumption on $l$.
Therefore, the line $\varpi'(l')$ is not contained in the branch quartic of $\varpi'$.
Hence, the curve $(\varpi')\inv(\varpi'(l'))$ on $S'$ has a component
other than the line $l$ which is mapped to the bitangent $\varpi'(l')$
isomorphically.
Since the branch curve is of degree four,
this can happen exactly when the line $\varpi'(l')=\varpi(l)$ is a bitangent
of the branch quartic.
Now we move the line along the branch quartic in a way that 
it is a regular tangent of the branch quartic.
The preimages of such a tangent under $\varpi$ and $\varpi'$ are still
anti-canonical curves on $S$ and $S'$ respectively,
and generically have respectively an ordinary double point over the tangent point
as their all singularity.
This means that the preimage to $S$ of the tangent is a minitwistor line
on $S$.
Therefore, the curve $\varpi\inv(\varpi(l))$ on $S$ is  a limit of 
minitwistor lines on $S$.
Since $S^*$ is closed and $\varpi\inv(\varpi(l))= H\cap S$, this means $H\in S^*$.
Hence $l^*\subset S^*$, as desired.
\proofend

\medskip
From the proof, it is not difficult to obtain the following

\begin{figure}
\includegraphics{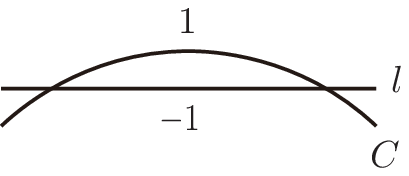}
\caption{
}
\label{fig:cycle1}
\end{figure}

\begin{proposition}
If $l$ is a line on a Segre surface $S\subset\CP_4$ 
which does not pass any singularity of $S$ as in 
Proposition \ref{prop:line}, 
for a generic hyperplane $H$ which contains $l$,
the section $S\cap H$ is as in Figure \ref{fig:cycle1},
where the numbers are self-intersection numbers in $S$.
\end{proposition}

\proof
We keep the notations in the previous proof.
The image $\varpi(l) = \varpi'(l')$ was a bitangent of the branch curve of
the double covering $\varpi':S'\to \CP_2$, and
the preimage $(\varpi')\inv(\varpi'(l'))$ is of
the form $l' + C'$ for some smooth rational curve $C'\subset S'$
which is mapped isomorphically to the bitangent by $\varpi'$.
Moreover the intersection $l'\cap C'$ consists of two points,
which are over the tangent points of the bitangent,
and both intersections are transversal.
Let $C$ be the image of $C'$ by the blow-down $S'\to S$.
Since $S|_H = l + C$, $C$ is a cubic curve and by rationality 
it is a rational normal curve in $H$.
Also $C$ is contained in the smooth locus of $S$,
and the intersection $C\cap l$ consists of two points
and both intersections are transversal.
For finishing a proof of the proposition,
since $l$ is a $(-1)$-curve, it remains to
see that $C^2 = 1$.
Since $C+l$ is a hyperplane section of $S$,
we have $(C+l)^2 = \deg S = 4$.
Moreover we have $C.\,l = 2$ as above. 
This means $C^2 = 1$, as desired.
\proofend

\medskip
Next we find another type of divisors at infinity.
For this purpose we use the double covering structure
over a smooth quadric obtained in 
Proposition \ref{prop:SS}.
So let $S$ be any one of  Segre surfaces listed in Table \ref{table1}.
As in the proof of Proposition \ref{prop:SS} 
we denote $\pi:\CP_4\to\CP_3$ for the generating projection
from a point that induces the double covering map
from $S$
to a smooth quadric $Q\subset\CP_3$.
We denote the generating point of $\pi$ and its dual hyperplane
respectively by
\begin{align}\label{wwd}
w\qandq w^*.
\end{align}
In the coordinates of the proof of Proposition \ref{prop:SS}, 
$w = (1,0,0,0,0)$.
The projection $\pi$ induces the dual inclusion 
\begin{align}\label{pid}
\pi^*:\CP_3^*\hookrightarrow \CP_4^*,
\end{align}
and we have $\pi^*(\CP_3^*) = w^*$.
In the following we often think subvarieties in $\CP_3^*$
as those in $\CP_4^*$ through this inclusion.

If $B\subset\CP_3$ is an irreducible curve and is not a straight line,
the dual variety $B^*\subset\CP_3^*$ is defined
in the same way to the case for the surface $S$
explained in Section \ref{s:mts},
and it is irreducible and 2-dimensional.
If $B\subset\CP_3$ is reducible and does not have a straight line as a component, we define the dual variety $B^*$
to be the union of the dual varieties of all its irreducible components.
If $B$ is reducible and have a straight line as a component, we define $B^*$ to be the union of
the dual varieties of all its irreducible components
which are not a straight line.
So $B^*$ is empty if all components of $B$ are lines.


\begin{proposition}\label{prop:dcQ}
Assume that a Segre surface $S$ is a finite
double cover over a smooth quadric surface $Q$
by 
a projection from a point of $\CP_4$
as in Proposition \ref{prop:SS}.
Let $w\in\CP_4$ and $w^*\subset\CP_4^*$ be as in \eqref{wwd}.
Then if $B\subset Q$ denotes the branch divisor
of the double cover,
the hyperplane section by $w^*$ 
of the dual variety $S^*$ satisfies
\begin{align}\label{hpsc1} 
  Q^* \cup B^*\subset S^*\cap w^*,
\end{align}
where $Q^*$ and $B^*$
are the dual varieties of $Q$ and $B$ respectively.
\end{proposition}

\proof 
The inclusion $Q^*\cup B^*\subset w^*$ is  obvious
because we are using the inclusion \eqref{pid},
so it is enough to show $Q^*\cup B^*\subset S^*$.

First, we show $B^*\subset S^*$.
From the above convention for the dual varieties,
it is enough to show that for any irreducible component $B_1$ of $B$
which is not a straight line, the inclusion $B_1^*\subset S^*$ holds.
In the following for simplicity we write $B$ for $B_1$.
Then since $B$ is not a line,
as a curve on $Q\simeq\qdr$, it is not a
curve of bidegree $(1,0)$ nor $(0,1)$.
Therefore, a smooth generic point $q$ of $B$
satisfies $T_qB\not\subset Q$.
For such a point $q\in B$, 
the intersection of $Q$ with a 2-plane $h\subset\CP_3$ containing $T_qB$
is a smooth $(1,1)$-curve on $Q$ unless $h = T_qQ$.
Moreover except for a finite number of such 2-planes, the intersection
$h\cap B$ consists of three points, and one of them is the tangent point $q$,
while the other two intersections are transversal.
If we put $H:=\pi\inv(h)$, then
the hyperplane section $S\cap H$ has a node at the point $(\pi|_S)\inv (q)$
because the curve $Q\cap h$ is tangent to $B$ at $q$,
and it has no other singularity
because $h$ intersects transversally at the other two intersection points with $B$.
This means $H\in W_0$, where as before $W_0$ is the Zariski-open subset of $S^*$ which parameterizes minitwistor lines
(see Definition \ref{def:Severi} and Proposition \ref{HNMT}).
Hence, for a generic $h\in B^*$, we have $\pi\inv(h)\in S^*$.
Since $B^*$ is closed, this implies  $B^*\subset S^*$.

%

It remains to show $Q^*\subset S^*$.
So let $h\subset\CP_3$ be a tangent plane to $Q$,
and $q$ the tangent point.
The section $Q\cap h$ consists of two lines
intersecting transversally at the point $q$.
Let $l_1$ and $l_2$ be these lines.
We choose the tangent plane $h$ in such a generic way that $q\not\in B$
and that the two lines $l_1$ and $l_2$ intersect
the branch divisor $B$ transversally.
Since $B$ is a cut of $Q$ by a quadric 
by Proposition \ref{prop:SS},
the intersections $B\cap l_1$ and $B\cap l_2$
consist of two points respectively.
From these the preimages 
\begin{align}\label{C1C2}
C_1:=(\pi|_S)\inv(l_1)\qandq C_2:=(\pi|_S)\inv(l_2)
\end{align}
are smooth rational curves on $S$, and they intersect
transversally at the two points $(\pi|_S)\inv(q)$.
We write $p_1$ and $p_2$ for these two points,
and put $H=\pi\inv(h)$. The latter is a hyperplane
spanned by the generating point $w$ of the 
projection $\pi$ and the tangent plane $h$,
and contains the tangent spaces $T_{p_1}S$ and $T_{p_2}S$
since $T_qQ = h$ and $q\not\in B$.
Obviously we have $H|_S = C_1 + C_2$,  
$C_1 \cap C_2 = \{p_1,p_2\}$.
Moving the tangent point $q\in Q$,
we obtain a 2-dimensional family of reducible curves belonging to 
the system $\big|\ms O_{\CP_4}(1)|_S\big| = \big| K_S\inv \big|$,
whose members have  two nodes
as their all singularities. 
In particular, $H\not\in W_0$.

In order to show $Q^*\subset S^*$ it is enough to
show that the curve $C_1+C_2$ is a limit of 
curves which belong to the Zariski-open subset $W_0$ of $S^*=W$. 
Since $l_1^2 = l_2^2 = 0$ on $S$,  from \eqref{C1C2},
we have $C_1^2 = C_2^2 = 0$ on $S$.
Hence by adjunction we have  $K_S.\,C_1 = K_S.\,C_2 = -2<0$.
From \cite[Prop.\,(2.11)]{Tan},
this means that any one of the two nodes
$p_1$ and $p_2$ of the curve $C_1+C_2$ can be 
smoothed out under a small displacement by moving the curve on $S$
while the other node is not smoothed out.
Clearly the curve
obtained as such a partial smoothing of the curve $C_1+C_2$
belongs to the subset $W_0$ of $S^*$.
This means that $C_1+C_2$ is a limit of curves which belong to $W_0$,
and so we obtain $Q^*\subset S^*$, as desired.
\proofend

\medskip
Note that from the proof, generic points of $B^*$ belong to $W_0$.
Therefore, $W_0$ is not entirely contained in the affine space
$\CP_4^*\minus w^*\,(\simeq\CC^4)$.
The following proposition is also obvious from the proof.

\begin{figure}
\includegraphics{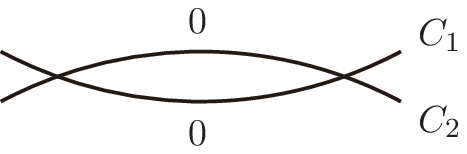}
\caption{
}
\label{fig:cycle2}
\end{figure}

\begin{proposition}
In the situation of the previous proposition,
for a generic hyperplane $H$ that belongs to the dual quadric
$Q^*$, the section $S\cap H$ is as in Figure \ref{fig:cycle2}.
\end{proposition}

Thus, we have found two kinds of divisors at infinity,
and both consist of two smooth rational curves
intersecting transversally at two points.
These two kinds of curves are
distinguished by 
the self-intersection numbers of the two components.

We next show that 
the dual variety $S^*$ intersects itself along 
the dual plane $l^*$ and the dual quadric $Q^*$.
More precisely, we show

\begin{proposition}\label{prop:si1}
Let $l\in S$ and $Q\subset\CP_3$ be
as in Propositions \ref{prop:line} and \ref{prop:dcQ}.
Then the dual variety $S^*$ has ordinary double points
along the dual plane $l^*$ and
the dual quadric $Q^*$.
\end{proposition}

\proof
We use results and notations in \cite{Tan}.
(In \cite{Tan} the letter $D$ is used instead of $C$.)
The proofs for $l^*$ and $Q^*$ are almost the same,
so for a generic hyperplane $H$ belonging to $l^*$ or $Q^*$
we write $H|_S = C =C_1+C_2$  and $C_1\cap C_2=\{p_1,p_2\}$.
By choosing $H\in l^*\cup Q^*$ in a sufficiently generic way,
we may suppose that $C\cap \Sing\,S = \emptyset$.
Put $N_C:=\ms O_C(C)$. This is an invertible sheaf
on the curve $C$, and 
there is a natural isomorphism
$$
T_H \big|\ms O_{\CP_4}(1)\big|
\simeq H^0(N_C)\simeq\CC^4.
$$
Moreover from the exact sequence 
\begin{align}\label{osc}
0 \lras \ms O_S \lras \ms O_S(C) \lras 
N_C \lras 0
\end{align}
and rationality of $S$, we have
$H^1(N_C)\simeq H^1\big(\ms O_S(C)\big)$.
Let $\tilde S\to S$ be the minimal resolution of all singularities of $S$
and $\tilde C$ the strict transform of $C$ into $\tilde S$.
Then $\tilde C$ is an anti-canonical curve on $\tilde S$, and it is nef and big.
Hence by Kodaira-Ramanujan vanishing,
we have $H^1\big(\ms O_{\tilde S}(\tilde C)\big)=0$.
So as $C\cap \Sing S = \emptyset$,
we have $H^1\big(\ms O_S(C)\big) = 0$, and hence from \eqref{osc}
we obtain $H^1(N_C) = 0$.

Next write $T^1_C:=\ms E\!xt^1_{\ms O_C}
(\Omega_C,\ms O_C)$, where $\Omega_C$ is
the sheaf of K\"ahler differentials on $C$.
Then from the standard exact sequence 
$0 \to \ms O_S(-C)\stackrel{d}{\to} \Omega_S|_C
\to \Omega_C\to 0$,
there is a natural homomorphism 
$N_C\to T^1_C$.
Since $C$ is normal crossing at $p_1$ and $p_2$, 
we have $T^1_C  \simeq \CC_{p_1}\oplus \CC_{p_2}$
and since $S$ is smooth at these points, the homomorphism $N_C\to T^1_C$ is surjective.
Let $N'_C$ be the kernel sheaf of this homomorphism,
so that we have an exact sequence
\begin{align}\label{ext1}
0 \lras N'_C \lras N_C \lras T^1_C\lras 0.
\end{align}
The Zariski-tangent space at the point $[C]$
of the locus in the linear system $|K_S\inv| \simeq \big|\ms O_{\CP_4}(1)\big|$
which corresponds to equisingular displacements of $C$
is given by $H^0(N'_C)$, 
and the obstruction for smoothness of this locus is in 
$H^1(N'_C)$.

Now in case $H\in l^*$, namely in case $l\subset H$,
as in the proof of Proposition \ref{prop:line}, 
the image $\varpi(l)$ is a bitangent of the branch quartic 
of the double cover $\varpi':S'\to\CP_2$
obtained by choosing a secant $\ol{pq}$ as in Lemma \ref{lemma:scnt},
and the two singularities $p_1$ and $p_2$ of the curve $C$
are mapped to the tangent points of the bitangent.
Evidently there are two ways of moving the bitangent to a regular tangent
by respecting the tangency at $p_1$ or $p_2$.
Accordingly, the reducible curve $C=C_1+C_2$ admits a displacement
in $S$ for which exactly any one of the singularities
$p_1$ and $p_2$ is smoothed out.
This means that in case $H\in l^*$
the natural map $H^0(N_C)\to H^0(T^1_C)\simeq\CC^2$ is surjective.
In case $H\in Q^*$, as in the last part of the proof of 
Proposition \ref{prop:dcQ},
the curve $C$ admits a displacement in $S$
which induces a smoothing for exactly any one of 
the singularities $p_1$ and $p_2$. 
Hence 
the natural map $H^0(N_C)\to H^0(T^1_C)\simeq\CC^2$
is again surjective.
Hence, regardless of $H\in l^*$ or
$H\in Q^*$, from \eqref{ext1} and $H^1(N_C)=0$, we obtain $H^1(N'_C)=0$.
We also obtain $h^0(N'_C) = h^0 (N_C) - h^0(T^1_C)= 4-2=2$.

Finally, we define the sheaves $N''_C(p_1)$ and $N''_C(p_2)$ on $C$ by the properties
\begin{align*}
N''_C(p_1)|_{C \minus \{p_1\}}
\simeq N_C|_{C \minus \{p_1\}}
\qandq
N''_C(p_1)|_{C \minus \{p_2\}}
\simeq N'_C|_{C \minus \{p_2\}},\\
N''_C(p_2)|_{C \minus \{p_2\}}
\simeq N_C|_{C \minus \{p_2\}}
\qandq
N''_C(p_2)|_{C \minus \{p_1\}}
\simeq N'_C|_{C \minus \{p_1\}}.
\end{align*}
Since $N_C$ and $N'_C$ are canonically isomorphic
over $C\minus\{p_1,p_2\}$, these two sheaves are 
well-defined.
Note that the these sheaves
are slightly larger than $N'_C$,
and if $\{i,j\} = \{1,2\}$,
we have exact sequences
\begin{align}\label{ext2}
0 \lras N'_C \lras N''_C(p_i) \lras \CC_{p_j}\lras 0.
\end{align}
As in \cite[Remark (1.7)]{Tan},
for $i=1,2$, the space $H^0\big(N''_C(p_i)\big)$ is 
the Zariski-tangent space at the point $[C]$
of the locus in $\big|\ms O_S(C)\big|$ which corresponds to
equisingular displacements of $C$, where equisingularity is imposed only for the node $p_i$.
Moreover the space $H^1\big(N''_C(p_i)\big)$ is the obstruction
space for such displacements.
We then have  natural inclusions
\begin{align}\label{ext3}
H^0(N'_C)\subset H^0\big(N''_C(p_i)\big)\subset H^0(N_C),\quad i=1,2,
\end{align}
and the codimensions of the two inclusions are both one.
Since the nodes $p_1$ and $p_2$ can be independently smoothed out
as above, we have 
$H^0\big(N''_C(p_1)\big)\neq H^0\big(N''_C(p_2)\big)$
as subspaces in $H^0(N_C)$.
Hence we obtain the transversality
\begin{align}\label{transv}
H^0\big(N''_C(p_1)\big)\cap H^0\big(N''_C(p_2)\big) = H^0(N'_C).
\end{align}
Since $H^1(N'_C) = 0$ as above, from \eqref{ext2},
we obtain $H^1\big(N''_C(p_i)\big) = 0$ for $i=1,2$.
Thus, the above two equisingular displacements of $C$ in $S$ are 
unobstructed, and each partially equisingular deformations
constitute a smooth threefold in 
the linear system $\big|\ms O_S(C)\big|$.
By the transversality \eqref{transv}, 
these two components intersect transversally
along the locus of equisingular displacements of
$C$, where this time equisingularity is imposed on both 
$p_1$ and $p_2$.
Evidently the last locus is identified with a neighborhood of 
the point $H$ in the dual 2-plane $l^*$ or
the dual quadric $Q^*$.
Therefore, the dual variety $S^*$ has 
ordinary double points along $l^*$ or $Q^*$.
\proofend

\medskip
When the quadric $\pi(S)\subset\CP_3$ is not smooth but
the cone over an irreducible conic as in Proposition \ref{prop:SC},
we write $\Lmd$ for an irreducible conic,
and $\Cone(\Lmd)\subset\CP_3$ for the cone over $\Lmd$.
Also, similarly to the notation in Proposition \ref{prop:dcQ}, we write $w$ for 
the generic point of the projection which induces
the double covering map  $S\to \Cone(\Lmd)$,
and $w^*\subset\CP_4^*$ for the dual hyperplane to $w$.
Then the following proposition can be shown in the same way to 
the assertion $B^*\subset S^*\cap w^*$ in Proposition \ref{prop:dcQ}.

\begin{proposition}\label{prop:dcC}
Assume that a Segre surface $S$ is a finite
double cover over $\Cone(\Lmd)$ by 
a projection from the point $w\in \CP_4$
as in Proposition \ref{prop:SC}.
Then if $B\subset \Cone(\Lmd)$ denotes the branch divisor
of the double cover, then
the hyperplane section by $w^*$ 
of the dual variety of $S$ satisfies
\begin{align}\label{hpsc2}
 B^*\subset S^*\cap w^*.
\end{align}
\end{proposition}


\section{Special hyperplane sections of the dual varieties}
\label{s:shps}
As is assumed so far, let $S\subset\CP_4$ be a Segre quartic
surface, and suppose that the Segre symbol of $S$
includes at least one `1'.
Then by Propositions \ref{prop:SS} and \ref{prop:SC},
$S$ is a double cover of a quadric surface
by a projection from a point $w\in\CP_4$.
By Propositions \ref{prop:dcQ} and \ref{prop:dcC},
the hyperplane section $S^*\cap w^*\subset\CP_4^*$ contains
the dual varieties $Q^*\cup B^*$ when
the quadric surface is smooth,
and $B^*$ when the quadric surface is the cone over a conic.
In this section, we first investigate the dual varieties of
irreducible quartic curves in $\CP_3$, and next 
by using them as well as the class formula obtained in Section \ref{ss:cf},
we determine these hyperplane sections in concrete forms.

\subsection{Dual varieties of irreducible quartic curves
in $\CP_3$}\label{ss:de}
We begin with the case where the curve is smooth.
\begin{proposition}\label{prop:e1}
If $B\subset\CP_3$ is a smooth complete
intersection of two quadrics,
then the dual variety $B^*\subset
\CP^*_3$ is an irreducible non-normal surface of degree $8$.
Moreover, $B^*$ intersects itself along
four irreducible conics in $\CP_3$
which are mutually disjoint.
\end{proposition}

\proof
The dual variety $B^*$ is irreducible since $B$ is irreducible.
For the remaining assertion,
the easiest proof would be to use normalized equations of the complete intersection as follows.
The Segre symbol for the smooth complete intersection
$B$ as in the proposition is $[1111]$,
and if $X_1,X_2,X_3,X_4$ are homogeneous coordinates on $\CP_3$,
the normalized quadratic equations of $B$ (see Section \ref{ss:segs}) are given by
\begin{align}\label{defq1}
\sum_{i=1}^4 \aaa_i X_i^2 = \sum_{i=1}^4 X_i^2 =0,
\end{align}
where the constants $\aaa_1,\aaa_2,\aaa_3$ and $\aaa_4$ 
are mutually distinct.
Then for any $i=1,2,3,4$, eliminating the variables $X_i$ from these two equations,
the elliptic curve $B$ is contained in the quadric
in $\CP_3$ defined by
\begin{align}\label{cone1}
\sum_{j\neq i}(\aaa_j-\aaa_i)X_j^2=0.
\end{align}
Since $\aaa_i\neq \aaa_j$ if $i\neq j$, 
this is the cone over an irreducible conic.
Thus, the elliptic curve $B$ is contained in
the four quadratic cones \eqref{cone1},
$1\le i\le 4$.
For each $i=1,2,3,4$, we write $e_i$
for the point of $\CP_3$ whose non-zero
entry in the homogeneous coordinates is $i$-th one only, and 
$\pi_i:\CP_3\to\CP_2$ for the generic projection 
from the point $e_i$.
Further we put
$
\Lmd_i:=\pi_i(B)
$
which are irreducible conics.
Of course, this is nothing but the conic defined by
the equation \eqref{cone1}.
The restriction of the projection $\pi_i$ 
to the elliptic curve $B$ is a degree-two morphism
onto the smooth rational curve $\Lmd_i$,
and therefore, by Hurwitz formula, it has exactly four branch points.
If $l$ is a  tangent line to the conic
$\Lmd_i$, and if the tangent point is not a branch point of
the double cover $B\to\Lmd_i$,
the preimage $\pi_i\inv(l)$ is a 2-plane which is tangent to the elliptic curve $B$ at the two points
over the tangent point.
Hence, we have $\Lmd_i^*\subset B^*$ for $i=1,2,3,4$.

The incidence variety
$I(B)=
\big\{(p,h)\in B\times\CP_3^*\set T_pB\subset h\big\}
$
is a $\CP_1$-bundle over $B$ and hence it is 
a smooth surface.
From the double fibration which is analogous to \eqref{diagram:fac1},
the dual variety $ B^*$ is obtained from
$I(B)$ by identifying 
two points $(p_1,h_1)$ and $(p_2,h_2)$ of $I(B)$ for which  the coincidence $h_1=h_2$ occurs;
namely exactly when $h_1=h_2$ is a bitangent plane of $B$
and $p_1$ and $p_2$ are its tangent points.
In particular, if $h$ is a bitangent plane to $B$ 
and 
$p_1$ and $p_2$ are tangent points,
then the point $h\in \CP_3^*$, which belongs to
$ B^*$ from the definition of the dual variety, is a non-normal point of $ B^*$,
and the natural projection 
$
I(B)\to  B^*
$
is the normalization of $ B^*$.
Therefore, the above tangent plane $\pi_i\inv(l)$
is a non-normal point of $ B^*$ for any tangent line $l$ of $\Lmd_i$.
Thus, the dual variety $B^*$
is non-normal at any points of the dual conics $\Lmd_i^*$, $1\le i\le 4$, considered as bitangent planes to $B$ in $\CP_3$
through the projections $\pi_i:\CP_3\to\CP_2$.

Next, we show that these four conics $\Lmd_i^*\subset\CP_3^*$ are mutually disjoint.
For simplicity of presentation, we only prove 
$\Lmd_1^*\cap \Lmd_2^*=\emptyset$.
Other cases can be shown in the same way.
If a plane $h\subset\CP_3$ belongs to $\Lmd_1^*\cap \Lmd_2^*$,
it passes the two points $e_1$ and $e_2$, so
$h$ is defined by an equation of the form $aX_3 + b X_4 = 0$,
$a,b\in\CC$.
This equation can be regarded as an equation of the tangent
lines $\pi_1(h)$ and $\pi_2(h)$.
On the other hand, equations of tangent lines of the conics
\eqref{cone1} which are of the form $a X_3 + b X_4 = 0$
have to be 
\begin{gather}\label{tan1}
(\aaa_3 - \aaa_1) \xi_3 X_3
+ (\aaa_4 - \aaa_1) \xi_4 X_4 = 0
\quad{\text{for $\Lmd_1$}},\\
(\aaa_3 - \aaa_2) \tilde {\xi}_3 X_3
+ (\aaa_4 - \aaa_2)  \tilde {\xi}_4 X_4 = 0
\quad{\text{for $\Lmd_2$,}}
\label{tan2}
\end{gather}
where $\xi_3,\xi_4,,\tilde{\xi}_3$ and $\tilde{ \xi}_4$ are solutions of
the equations
\begin{align}\label{tan3}
(\aaa_3 - \aaa_1) \xi_3^2 + (\aaa_4 - \aaa_1) \xi_4^2 
=
(\aaa_3 - \aaa_2) \tilde {\xi}_3^2 + (\aaa_4 - \aaa_2) \tilde {\xi}_4^2 
=0.
\end{align}
The points $(X_2:X_3:X_4) = (0:\xi_3:\xi_4)$
and $(X_1:X_3:X_4) = (0:\tilde{\xi}_3:\tilde{\xi}_4)$
are tangent points to $\Lmd_1$ and $\Lmd_2$ respectively.
These are concretely written down from \eqref{tan3},
and if we substitute them to \eqref{tan1} and \eqref{tan2},
we find that the two equations \eqref{tan1} and \eqref{tan2}
cannot be proportional, by using $(\aaa_1-\aaa_2)(\aaa_3 - \aaa_4)\neq 0$.
Hence, we obtain $\Lmd_1^*\cap \Lmd_2^*=\emptyset$.

Finally, we show that the dual variety $ B^*$ is of degree $8$.
We fix any $i$, $1\le i\le 4$,
and consider the above projection $\pi_i:B\to\Lmd_i$.
Let $\lmd_1,\lmd_2,\lmd_3$ and $\lmd_4$ be
the branch points of this double cover $\pi_i$.
For each $j=1,2,3,4$ we write $\lmd_j^*$ for the set of lines on $\CP_2$
which pass the branch point $\lmd_j$.
These are lines in the dual plane $\CP_2^*$ in which 
the dual conic $\Lmd_i^*$ lies.
If $l\in \lmd_j^*$, then the plane $\pi_i\inv(l)$ passes the ramification point on $B$
over the branch point $\lmd_j$, and is tangent to $B$ at
the ramification point because $\pi_i\inv(l)$ contains
the generating line of the cone which passes the ramification point.
Hence, $\pi_i\inv(l)\in  B^*$.
Thus, we have the inclusion $\lmd_j^*\subset  B^*$ for $j=1,2,3,4$
under the dual embedding $\pi_i^*:\CP_2^*\hookrightarrow \CP_3^*$.

We still fix the index $i$,
and let $e_i^*$ be the 2-plane in the  space $\CP_3^*$ (which is dual
to $\CP_3$ in which $B$ is embedded) formed by 
planes in $\CP_3$ which pass the point $e_i$.
What we have seen so far is that we have the inclusion
\begin{align}\label{hps1}
\Lmd_i^*\cup(\lmd_1^*\cup\lmd_2^*\cup\lmd_3^*\cup
\lmd_4^*) \subset  B^*\cap e_i^*
\end{align}
for the hyperplane section of $ B^*$.
All components in LHS lie on the dual plane to $\CP_2$
in which the conic $\Lmd_i$ lies, and they are considered
as subvarieties in $\CP_3^*$ through the dual inclusion $\pi_i^*$.
Now from the above description of
the dual conic $\Lmd_i^*$ as a non-normal locus
of $ B^*$, the conic $\Lmd_i^*$ is contained in $ B^*$
with multiplicity two.
On the other hand, since generic members of 
the lines $\lmd_1^*,\dots,\lmd_4^*$ are tangent to $B$ 
at precisely one point (which is the ramification point), the variety $ B^*$ 
is smooth at generic points on these lines.
In particular, the lines $\lmd_j^*$, $j=1,2,3,4$,
are contained in the hyperplane section $B^*\cap e_i^*$ with multiplicity one.
Hence, from \eqref{hps1} we have $$\deg B^*\ge 2\cdot 2 + (1+1+1+1) = 8.$$
To show that this is an equality,
it is enough to show that the inclusion
\eqref{hps1} is full.
This is equivalent to showing that the preimage of
any line on $\CP_2$ which is not tangent
to the conic $\Lmd_i$ and which does not 
pass the branch points $\lmd_j$, $1\le j\le 4$,
does not belong to $ B^*$.
But for a such line $l$, the intersection $l\cap\Lmd_i$
consists of two points, and these points are
not the branch points $\lmd_j$.
Hence, the intersection $\pi_i\inv(l)\cap B$
consists of distinct four points and all the intersections are transversal.
This means $\pi_i\inv(l)\not\in B^*$.
\proofend

\medskip
Next we consider the case where the curve
$B$ is a nodal rational curve.

\begin{proposition}\label{prop:e2}
If a curve $B\subset\CP_3$ is a complete
intersection of two quadric and is 
a 1-nodal rational curve,
then the dual variety $ B^*\subset \CP_3^*$ is an irreducible non-normal surface 
of degree $6$. Further the variety $ B^*$ intersects
itself along two smooth conics
in $\CP_3^*$.
\end{proposition}

\proof
Again, the irreducibility of $B^*$ follows from that of $B$.
Since $B$ is a 1-nodal rational curve,
the Segre symbol for $B$ is $[112]$ (see Section \ref{ss:dcS}).
As explained in Section \ref{ss:segs}, this means that
equations of $B$ can be supposed to be of the form
\begin{align}\label{B1}
\aaa_1 X_1^2 +\aaa_2 X_2^2 +2\aaa_3 X_3X_4 + X_4^2
=
X_1^2 + X_2^2 + 2 X_3X_4 = 0,
\end{align}
where $\aaa_1,\aaa_2$ and $\aaa_3$ are mutually distinct.
By eliminating the variable $X_1$ or $X_2$
from these two equations,
we obtain equations of two conics, or equivalently
those of the cones over the conics.
Again, these conics are readily seen to be irreducible.
Let $\Lmd_1$ and $\Lmd_2$ be the two conics obtained this way.
In the same way to the proof of the previous proposition, for $i=1,2$,
by the generic projection $\pi_i$ from the point $e_i\in\CP_3$,
the curve $B$ has a structure of double covering
over $\Lmd_i$, and 
the branch consists of three points,
exactly one of which is a double point.
Then by the same reason to the previous proof,
we have the inclusion
\begin{align}\label{12B}
\Lmd_1^*\cup\Lmd_2^*\subset B^*,
\end{align}
and the dual variety $B^*$ intersects itself along these two conics.

As in the previous proof, let $e_1^*$ be the 2-plane in the dual space $\CP_3^*$ formed by planes passing
the vertex $e_1$ of the cone over $\Lmd_1$.
Then from \eqref{12B} we have
\begin{align}\label{incl2}
\Lmd_1^*\subset B^*\cap e_1^*
\end{align}
and the conic $\Lmd_1^*$ is contained in RHS with multiplicity two.
Let $\lmd_1,\lmd_2$ and $\lmd_3$ be the branch points of the 
projection $B\to \Lmd_1$, and suppose that $\lmd_3$ is the double point.
So the node of $B$ is over the point $\lmd_3$.
Then by the same reason to the previous proof, we have the inclusion
$$
\lmd_1^*\cup\lmd_2^*\subset  B^*\cap e_1^*.
$$
Combined with \eqref{incl2}, we obtain 
\begin{align}\label{incl3}
\Lmd_1^*\cup (\lmd_1^*\cup\lmd_2^*)\subset  B^*\cap e_1^*.
\end{align}
We show that this inclusion is full.
Take any line $l\subset\CP_2$ which does not belong to 
LHS, and put $h = \pi_1\inv(l)$.
It is easy to see $h\not\in  B^*$ when  $l\not\in\lmd_3^*$.
So suppose $\lmd_3\in l$ and $\lmd_1,\lmd_2\not\in l$.
Then the line $l$ and the conic $\Lmd_1$ intersect precisely at
two points, and one of them is $\lmd_3$.
Let $\lmd$ be another intersection point.
Of course $\lmd\neq\lmd_1,\lmd_2$.
The intersection $\pi_1\inv(l)\cap B$ consists of three points,
two of which are over $\lmd$,
and the remaining point is the node of $B$.
The generating line $\pi_1\inv(\lmd)$ of the cone $\pi_1\inv(\Lmd_1)$
intersects
$B$ transversally at the two points.
Another generating line $\pi_1\inv(\lmd_3)$
does not intersect $B$ transversally,
but intersects each of the two branches of $B$ at the node transversally.
Because transversality is an open condition,
if we move the plane $h$ in a way that it avoids
the node of $B$, the moved plane still intersects $B$
transversally and this time at four distinct points.
This implies that the original plane $h$ does not belong to $ B^*$.
Thus, the inclusion \eqref{incl3} is an equality.
Hence, we obtain 
\begin{align*}
\deg  B^* &= \deg  (B^*\cap e_1^*) =
2\deg \Lmd_1^* + (1+1) = 2\cdot 2 + 2 = 6,
\end{align*}
as desired.\proofend

\begin{proposition}\label{prop:e3}
If a curve $B\subset\CP_3$ is a complete
intersection of two quadrics and is 
a rational curve having one cusp as its 
all singularities, then
the dual variety $ B^*\subset \CP_3^*$ is an irreducible non-normal surface 
of degree $5$. Further the variety $ B^*$ intersects itself along
an irreducible conic in $\CP_3^*$.
\end{proposition}

\proof 
Since this is almost parallel to 
the proofs of Propositions \ref{prop:e1} and \ref{prop:e2},
we just give an outline.
The Segre symbol of $B$ is $[113]$,
and from this we obtain normalized equations of $B$.
Eliminating one particular variable from the equations,
we obtain that $B$ is a complete intersection of 
the cone over an irreducible conic with another quadric.
If $\pi:B\to\Lmd$ denotes the projection to the conic,
this is two-to-one and has exactly two branch points.
One of them is a triple point and another one is a simple point.
The cusp of $B$ is over the triple point.
If $\lmd\in\Lmd$ is the simple branch point and $e_1$ is the vertex
of the cone,
by the same reason to the proofs of Propositions \ref{prop:e1} and \ref{prop:e2},
we have
$$
 \Lmd^* \cup \lmd^* \subset B^*\cap {e_1^*}, 
$$
and $B^*$ intersects itself along the dual conic $\Lmd^*$.
Further by the same argument as in the final part
of the proof of Propositions \ref{prop:e2},
this inclusion is shown to be full.
From these we obtain 
$$
\deg B^* = 2\cdot 2 + 1 = 5.
$$
Hence, we obtain the desired conclusion.\proofend

\subsection{Structure of special hyperplane sections of the dual varieties}
\label{ss:shs}
We return to the problem of describing  hyperplane sections 
$S^*\cap w^*\subset\CP_4^*$, where $w\in\CP_4$ is the center of the projection 
which induces 
 a double covering structure on a Segre quartic surface
as in Propositions \ref{prop:dcQ} and \ref{prop:dcC}.
As a first step we show

\begin{lemma}\label{lemma:intm}
Let  $b_1,\dots,b_l$ be all singularities of 
the branch curve $B$.
When the image $\pi(S)\subset\CP_3$ is a smooth quadric
as in Proposition \ref{prop:dcQ},
we have an inclusion
\begin{align}\label{incls1}
S^*\cap w^*\subset Q^*\cup B^*\cup (b_1^*\cup
\dots\cup b_l^*).
\end{align}
When the image $\pi(S)\subset\CP_3$ is the quadratic cone
as in Proposition \ref{prop:dcC}, 
if $v$ denotes the vertex of the cone,
we have an inclusion
\begin{align}\label{incls2}
S^*\cap w^*\subset v^*\cup B^*\cup (b_1^*\cup
\dots\cup b_l^*).
\end{align}
\end{lemma}

\proof
We only give a proof for  \eqref{incls1}
since \eqref{incls2} can be shown in 
almost the same way.
To prove  \eqref{incls1}, it is enough to see that there is no
component of the hyperplane section $ S^*\cap w^*$  other than $Q^*, B^*$ nor $b_1^*,\dots,b_l^*$.
To see this, let $H\subset\CP_4$ 
be a hyperplane belonging to 
the hyperplane $w^*$.
This is equivalent to $w\in H$, 
and we can write $H=\pi\inv(h)$ 
where $h\subset\CP_3$ is a 2-plane.
Suppose $H\not\in Q^*\cup B^*\cup (b_1^*\cup
\dots\cup b_l^*)$.
This means $h\not\in Q^*\cup B^*
\cup (b_1^*\cup
\dots\cup b_l^*)$.
Then the intersection $h\cap Q$ is a smooth
$(1,1)$-curve on $Q\simeq\qdr$ as $h\not\in Q^*$,
and it intersects
the branch curve $B$ transversally at four points
as $h\not\in B^*\cup (b_1^*\cup
\dots\cup b_l^*)$.
Hence, the restriction of $\pi|_S:S\to Q$
to the section $Q\cap h$ is a double cover over a 
smooth rational curve branched at four points,
and all the branch points are of multiplicity one.
This means that the cut of $S$ 
by the hyperplane $H=\pi\inv(h)$ is a smooth elliptic curve. Hence, $H\not\in  S^*$.
So we have obtained the inclusion \eqref{incls1}.
\proofend

\medskip
When the quadric $\pi(S)\subset\CP_3$ is smooth,
we have the following conclusion.
\begin{proposition}\label{prop:swqb}
If a Segre surface $S\subset\CP_4$ is a double cover over
a smooth quadric $Q$ by the projection from a point $w$
as in Proposition \ref{prop:dcQ}, we have the coincidence
\begin{align}\label{swqb}
S^*\cap w^*=Q^*\cup B^*,
\end{align}
where we regard the dual varieties $B^*$
and $Q^*$ as included in the hyperplane $w^*\subset\CP_4^*$.
\end{proposition}

\proof
The case where the branch curve $B$ is smooth is immediate
from \eqref{incls1}, and we need to show
the case where $B$ is singular.
Since we know $S^*\cap w^*\supset Q^*\cup B^*$
from Proposition \ref{prop:dcQ},
taking the multiplicity of $Q^*$ from 
Proposition \ref{prop:si1} into account, 
we have an inequality 
$
S^*|_{w^*} \ge 2Q^* + B^*
$
as  divisors on $w^*=\CP_3$.
Hence, we obtain
\begin{align}\label{S4B}
\deg S^* \ge  4 + \deg B^*.
\end{align}
Obviously, the equality holds if and only if 
the intersection $S^*\cap {w^*}$ has no component besides
$Q^*$ and $B^*$.
So it is enough to show that the equality
holds in \eqref{S4B}.

From Table \ref{table1} there are 9 symbols for which the Segre surfaces
are double covers over a smooth quadric and 
the branch curve $B$ is singular.
If the Segre symbol of $S$ is $[1112]$, 
the elliptic fibration $\tilde f':\tilde S'\to\CP_1$
in Corollary
\ref{cor:elf} has exactly one singular fiber
which is not of type I$_1$, and it is 
of type I$_2$. Hence, by
Theorem \ref{thm:cf}, we obtain $\deg S^*= 12 - 2 = 10$.
On the other hand, from Proposition \ref{prop:SS}
the branch curve $B$ is a 1-nodal rational curve
and hence by  Proposition \ref{prop:e2}
we have $\deg B^* = 6$.
Hence, the equality holds in \eqref{S4B}.
If the Segre symbol for $S$ is $[113]$,
a single singular fiber of type III (i.e.\,two rational curves touching at a single point)
appears instead of the above singular fiber of type I$_2$,
and from Theorem \ref{thm:cf} we obtain $\deg S^* = 12 - 3=9$.
On the other hand, from Proposition \ref{prop:SS}
the branch curve $B$ is a 1-cuspidal rational curve
and hence by Proposition \ref{prop:e3}
we have $\deg B^* = 5$.
Therefore, again the equality holds in \eqref{S4B}.

The branch curve $B$ for the case where the Segre symbol is $[111(11)]$ consists of two irreducible conics, and hence we have $\deg B^* = 2+2=4$.
On the other hand, since the two conics intersect transversally at two points,
$S$ has precisely two singularities, and both of them are type A$_1$.
Hence, from Theorem \ref{thm:cf}
we have $\deg S^* = 12-2-2=8$.
Hence the equality in \eqref{S4B} still holds.
The same proof works for Segre surfaces
with symbols $[12(11)],[1(11)(11)],[11(12)]$
and $[1(13)]$.
If the symbol of $S$ is $[122]$, the branch curve
$B$ consists of a rational normal curve and
a line in $\CP_3$.
The dual variety of a rational normal curve in $\CP_n$
is $(n+1)$ (see \cite[Example 10.3]{Tev}), and hence
neglecting the line component as before, we have $\deg B^* = 4$.
On the other hand the elliptic fibration
in Corollary \ref{cor:elf} has exactly two 
singular fibers which are not of type I$_1$,
and they are of type I$_2$.
Hence, from Theorem \ref{thm:cf},
we have $\deg S^* = 12 - 2-2=8$.
Hence the equality again holds in \eqref{S4B}.
The same argument works for the remaining one
whose symbol is $[14]$.
\proofend

\medskip
Next, we consider the case where the base space
of the double cover is a cone as in 
Proposition \ref{prop:dcC}.
From Table \ref{table2}, there are exactly 9 kinds of Segre surfaces
which have a structure of double covering over the quadratic cone, and exactly one of them, 
namely $[(12)(11)]$, has two such structures.
Among these 10 kinds,
the branch curves do not pass the vertex of the cone
for exactly 5 kinds.
For all of these, the same argument as
in the previous proposition works by using 
an inequality $\deg S^*\ge \deg B^*$ from Proposition 
\ref{prop:dcC} instead of \eqref{S4B},
and we obtain the following

\begin{proposition}\label{prop:swqb2}
If a Segre surface $S$ is a double cover over
a quadratic cone by the projection from a point $w$
as in Proposition \ref{prop:dcC},
and if the branch curve $B$ does not pass
the vertex of the cone, we have the coincidence
\begin{align}\label{swqb2}
S^*\cap w^*= B^*,.
\end{align}
where we regard the dual variety $B^*$ as included in the hyperplane $w^*\subset\CP_4^*$.
\end{proposition}

In particular we have $v^*\not\subset S^*$ for 
the 2-plane $v^*$ which is dual to the vertex $v$ of 
the quadratic cone.
This means that if $H\subset\CP_4$ is a generic hyperplane 
which passes the two A$_1$-singularities 
that are over the vertex $v$ of the cone,
any small displacement of $H$ in $\CP_4$ does not
cut out a minitwistor line from $S$.

Next, in order to treat the case where 
the branch curve $B$ passes the vertex $v$ of the cone, 
we show the following

\begin{proposition}\label{prop:odp}
Suppose that a Segre quartic surface $S\subset\CP_4$ 
has an ordinary double point,
and let $H\subset\CP_4$ be a hyperplane
such that the intersection $S\cap H$ has a node
at the ordinary double point as its all singularity.
Then $H\not\in S^*$.
\end{proposition}

\proof
Let $p_i$ be an ordinary double point of $S$
and $H\subset \CP_4$ a hyperplane through $p_i$,
and suppose that the section $C:=S\cap H$
has a node at $p_i$ as its all singularity.
Let $\mu:\tilde S\to S$ be the blowing-up at $p_i$.
This resolves the singularity of $S$ and $C$.
Let $\tilde C$ be the strict transform of $C$ into $S$
and $E$ the exceptional curve.
Since $C$ is a 1-nodal rational curve,
$\tilde C$ is a smooth rational curve, and
we have  $\mu^*C = \tilde C+E$.
Further the curve $\tilde C$ intersects $E$ transversally at exactly two points.
So the reducible curve $\tilde C + E$ has two nodes as its all singularities.
If $H\in S^*$ would hold, then the curve $C$
on $S$ can be deformed into 
a minitwistor line under a small displacement in $S$.
Then by lifting it to $\tilde S$, 
we would obtain a small displacement of the curve $\tilde C + E$
whose resulting curve is 1-nodal irreducible and does not intersect $E$.
To show that $\tilde C + E$ 
does not admit such a small displacement,
we write $\ms I_{\tilde C + E}\subset\ms O_{\tilde S}$
for the ideal sheaf of the curve $\tilde C + E$,
and consider an exact sequence
\begin{align}\label{}
0 \lras \ms I_{\tilde C + E}/\ms I_{\tilde C + E}^2
\stackrel{d}{\lras} \Omega_{\tilde S}|_{\tilde C + E} 
\lras \Omega_{\tilde C + E} \lras 0.
\end{align}
Writing $T_X$ for the tangent sheaf of a variety $X$,
$T^1_X$ for the  Ext-sheaf $\ms E\!xt^1(\Omega_X,\ms O_X)$,
and $\mathbb T^1_X$ for the vector space ${\rm{Ext}}^1(\Omega_X,\ms O_X)$,
this induces the commutative diagram
\[
  \begin{CD}
     0 @>>> H^0\big(T_{\tilde C + E}\big) @>>> 
     H^0\big(T_{\tilde S}|_{\tilde C + E}\big) @>>> 
     H^0\big( \ms O_{\tilde C + E}(\tilde C + E)\big) @>{\aaa}>> 
     \mathbb T^1_{\tilde C + E}
 \\
&& @VVV @VVV  @VVV @VVV \\
     0   @>>>  T_{\tilde C + E} @>>> T_{\tilde S}|_{\tilde C + E} @>>> 
     \ms O_{\tilde C + E}(\tilde C + E) @>{\bbb}>> T^1_{\tilde C + E} 
  \end{CD}
\]
where each row is exact.
Let $p$ and $q$ be the singularities of the curve $\tilde C+E$.
Namely $\tilde C\cap E = \{p,q\}$.
Then since these are nodes of $\tilde C+E$
and both $\tilde C$ and $E$ are rational curves intersecting at two points, 
we have 
\begin{align}\label{T1iso}
T^1_{\tilde C+E}\simeq \CC_p \oplus \CC_q
\qandq
\mathbb T^1_{\tilde C+E} \simeq H^0(\CC_p \oplus \CC_q
) \simeq \CC^2.
\end{align}
Moreover the map $\bbb$ in the second row
is the map which just assigns the restrictions of the germs at the points $p$ and $q$.
Now since $(\tilde C+E).\, E=2-2=0$, 
the invertible sheaf $\ms O_{\tilde C + E}(\tilde C + E)$ is trivial on 
the component $E$.
Hence from the commutativity of the above diagram,
there is a natural identification between
the two skyscraper sheaves $\CC_p$ and $\CC_q$,
as well as the two components of $\mathbb T^1_{\tilde C + E}
\simeq H^0(\CC_p \oplus \CC_q) $.
%

The vector space 
$H^q\big( \ms O_{\tilde C + E}(\tilde C + E)\big)$ is 
the space of first-order displacements of $\tilde C+E$ in $\tilde S$
when $q=0$, and is the obstruction space for such displacements when $q=1$.
Moreover, from Kodaira-Ramanujan vanishing theorem
as in the proof of Proposition \ref{prop:si1},
we have $H^1\big( \ms O_{\tilde S}(\tilde C + E)\big)=0$
and this easily means
$H^1\big( \ms O_{\tilde C + E}(\tilde C + E)\big)=0$.
From these, it follows that 
a neighborhood of the origin in the space $H^1\big( \ms O_{\tilde C + E}(\tilde C + E)\big)$
is identified with a neighborhood of 
the Kuranishi family for displacements of $\tilde C+E$ in $\tilde S$.
From the above natural identifications over the singularities $p$ and $q$,
this means that  under a displacement in $\tilde S$, one of the two singularities $p$ and $q$ 
is also smoothed out
iff the other singularity is smoothed out.
Therefore, there does not exist a displacement of $\tilde C+E$ in $\tilde S$
which gives a smoothing of exactly one of $p$ and $q$.
Hence $H\not\in S^*$.
\proofend

\begin{proposition}\label{prop:swqb3}
If a Segre surface $S$ is a double cover over
a quadratic cone by the projection from a point $w$
as in Proposition \ref{prop:dcC},
and if the branch curve $B$ passes
the vertex $v$ of the cone, then we have 
\begin{align}\label{swqb3}
S^*\cap w^*= B^*\cup v^*,
\end{align}
where we regard the dual plane $v^*$ as included in $w^*$ in a natural way.
\end{proposition}

\proof
We first show that when $v\in B$,
one always has a discrepancy 
$\deg S^* > \deg B^*$.
As in Table \ref{table2} the situation $v\in B$ happens exactly for the following five kinds of Segre surfaces:
\begin{align}\label{vinB}
[(12)11],\,\,[1(13)],\,\,[(14)],\,\,[(12)2],\,\,[(11)(12)].
\end{align}
Here, for the final surface,
`1' is chosen from the ingredient of the first parenthesis.
When the symbol of a Segre surface
$S$ is $[(12)11]$, we have $\deg S^* = 8$.
On the other hand, since the branch curve $B$ on the cone
is a 1-nodal rational curve, we have $\deg B^* = 6$ 
by Proposition \ref{prop:e2}. So $\deg S^*>\deg B^*$ holds.
When the symbol of a Segre surface
$S$ is $[1(13)]$,  we have
$\deg S^* = 6$.
On the other hand, by Proposition \ref{prop:SC},
the curve $B$ is a cuspidal rational curve,
so we have $\deg B^* = 5$ from Proposition \ref{prop:e3}. 
Thus, again $\deg S^* > \deg B^*$ holds.
In the same way, when the Segre symbol is $[(14)]$
we have $\deg S^* = 5$ and $\deg B^* = 4$.
When the Segre symbol is $[(12)2]$
we have $\deg S^* = 6$ and $\deg B^* = 4$.
Finally, when  the Segre symbol is $[(11)(12)]$
and if `1' is chosen from the ingredient of the first parenthesis,
we have $\deg S^* = 4$ and $\deg B^* = 2$.
Hence, we have obtained the discrepancy $\deg S^*>\deg B^*$ for
all these five kinds of surfaces.

The discrepancy $\deg S^*> \deg B^*$ means that 
the inclusion $B^*\subset S^*\cap w^*$ obtained in 
Proposition \ref{prop:dcC} is not full.
Namely the intersection $S^*\cap {w^*}$ has at least one irreducible component other than $B^*$.
By \eqref{incls2} in Lemma \ref{lemma:intm}, 
that component is either the dual plane $v^*$ 
to the vertex $v$ of the cone or the dual plane
$b_i^*$ for some singularity $b_i$ of the curve $B$.
Among these five kinds of surfaces,
the branch curves $B$
for the three cases $[(12)11]$, $[1(13)]$ and $[(14)]$ have 
a singularity only at the vertex $v$ of the cone
because the singularity of $B$ is unique and is exactly the vertex.
Therefore, for these three kinds of surfaces,
 the `extra' component of $S^*\cap {w^*}$
has to be the dual plane $v^*$.
The branch curves for the remaining two kinds of surfaces 
have singularity which is not the vertex of the cone,
but all of them are ordinary nodes.
So by Proposition \ref{prop:odp}, 
the 2-planes which are dual to these nodes are not included in $S^*$.
Therefore, again the extra component of $S^*\cap w^*$ 
has to be the 2-plane $v^*$.
\proofend

\medskip
Among the five kinds of the Segre surfaces as in \eqref{vinB},
for the two kinds $[1(13)]$ and $[(14)]$,
as seen in the above proof, we have $\deg S^* - \deg B^*  = 1$.
Hence, the multiplicity of the plane $v^*$in the section $S^*\cap {w^*}$ is one. 
Thus, for these two kinds of surfaces, as divisors on $w^*=\CP_3$ we have
\begin{align}\label{Sw3}
S^*|_{w^*} = B^* + v^*.
\end{align}
For the rest of the surfaces, the ones whose symbols are
$[(12)11], [(12)2]$ and $[(11)(12)]$,
 we have
$\deg S^* - \deg B^*  = 2$, and therefore we obtain 
\begin{align}\label{2v}
S^*|_{w^*} = B^* + 2v^*.
\end{align}
Thus, we have determined all hyperplane sections
of the form $S^*\cap w^*$ in explicit forms.

\section{Examples of transitions and a concluding remark}
\subsection{typical transitions between Segre quartic surfaces}
\label{ss:dgn}
Next, we discuss adjacent relations (or degeneration relations) between
some Segre quartic surfaces.
First, the transitions $[1111]\to [112]\to [13]$ 
for complete intersections of two quadrics in $\CP_3$ provide
a standard degeneration of a smooth elliptic curve into cuspidal rational curve via a nodal rational curve, and accordingly, 
through the double cover,
the transitions $[11111]\to[1112]\to[113]$
for Segre surfaces
provide modest degenerations of the surfaces.
From Table \ref{table1}, in this degeneration, 
the classes of the Segre surfaces 
decrease as $12\to 10\to 9$,
while those of the branch curves decrease as
$8\to 6\to 5$ as obtained in 
Section \ref{ss:de}.
Thus, all decreases of the classes of Segre surfaces
 come from those of the branch curves.

For another typical degenerations of 
complete intersections of two quadrics in $\CP_3$,
we take 
\begin{align}\label{dege1}
[112]\to [11(11)]\to [2(11)]\to [(11)(11)].
\end{align}
These are $1$, $2$, $3$ and $4$-nodal curves respectively,
and the numbers of irreducible components of the curves are
1,2,3 and $4$ respectively.
The first one is a $1$-nodal rational curve,
and the second one consists of two conics
(i.e.\,$(1,1)$-curves in $Q\simeq\qdr$) intersecting
transversally at two points.
The third one consists of two lines
and one conic, and they form a `triangle' of smooth
rational curves.
The final one consists of four lines,
forming a `square' of rational curves.
Using that the dual variety of a smooth conic in $\CP_3$ is 
the cone over the conic and in particular quadratic,
we obtain that the degrees of the dual varieties of these
quartic curves are respectively 6, $2+2 = 4$, 2 
and $0$.
Adding a single `1` to each symbol in the series
\eqref{dege1},
we obtain a series of Segre surfaces whose symbols are
\begin{align}\label{tr1}
[1112]\to [111(11)]\to [12(11)]\to [1(11)(11)].
\end{align}
From the above description of the quartic curves,
these Segre surfaces
have one, two, three and four ${{\rm A}}_1$-singularities
respectively, and have no other singularities.
Hence by Theorem \ref{thm:cf},
the classes of these surfaces are $10,8,6,4$
respectively.
Thus, again the decreases of the classes
of the surfaces exactly come from those of the branch quartic curves.

For one more  interesting example of a degeneration,
we take the transition $[22]\to [4]$ for
complete intersections in $\CP_3$.
Both curves consist of one rational normal curve
 and one line in $\CP_3$,
 but
 for the symbol $[22]$ the two components intersect transversally at two points,
while for the symbol $[4]$ the two components are tangent
at one point.
If we consider the transition 
$[122]\to [14]$ for Segre surfaces, 
the singularities are two A$_1$ for the former and 
a single A$_3$ for the latter,
Hence, from Theorem \ref{thm:cf}, for both kinds of surfaces, we have
$\deg S^* = 12 -  4 = 8$ (see Table \ref{table1}).
This can be also obtained from Proposition \ref{prop:swqb}
if we note that the class of a rational normal curve
in $\CP_n$ is $(n+1)$ as mentioned before.
Thus, no decrease occurs for the classes in this degeneration.

Next, we consider the transition $[(11)3]\to
[(13)1]$ for Segre surfaces.
Both of these kinds of Segre surfaces have a structure of 
double covering over a quadratic cone as in Proposition \ref{prop:dcC}.
The branch curves have $[13]$ as the Segre symbol for both,
and this is a cuspidal rational curves.
The vertex $v$ of the cone
is on the branch curve only for the latter kind of surfaces.
So certainly, the latter kind of surfaces is  obtained 
from the former kind of surfaces as a degeneration.
The former kind of surface has two A$_1$-singularities and
one A$_2$-singularity, and hence the class is  $12-2-2-3=5$.
This can be also obtained from Propositions \ref{prop:swqb2} and \ref{prop:e3}.
On the other hand, the latter kind of surfaces
have exactly one singularity which is of type D$_4$. 
Hence, we have $\deg S^*= 12-6=6$.
Thus, the class increases under this degeneration,
and by Proposition \ref{prop:swqb3} it comes from the component $v^*$.
In other words, in the case of $[(13)1]$,
a hyperplane section of $S$ which passes the D$_4$-singularity
admits a displacement which avoids the singularity 
and which gives a 1-nodal curve which is really 
a minitwistor line.
This is in high contrast with the situation mentioned in Remark \ref{rmk:mtl}.

As a final example of a degeneration,
we consider the transition $[(12)2] \to [(14)]$.
Segre surfaces of these two kinds have again a double covering
structure over a quadratic cone, and 
the components of the branch curves on the cone are one generating
line and one rational normal curve.
Moreover, both of the two components pass the vertex of the cone.
The difference for these two branch curves is that, for $[(12)2]$ 
the two components intersect transversally
at two points, while for $[(14)]$
they are tangent at the vertex and have no other intersection.
From Theorem \ref{thm:cf} we have
$\deg S^* = 6$ in the case  $[(12)2]$, and
$\deg S^* = 5$ in the case  $[(14)]$.
As we have already obtained in \eqref{2v} and \eqref{Sw3} we have
$S^*\cap w^* = B^*+ 2v^*$ for $[(12)2]$, and
$S^*\cap w^* = B^*+ v^*$ for $[(14)]$.
Thus, the difference of the classes of the surfaces comes from 
the multiplicity of the dual plane $v^*$.

\subsection{Concluding remarks}\label{ss:cr}
We end this article by giving two remarks.
The first one is about a relation between null surfaces in 
Einstein-Weyl space associated to the Segre surfaces
and their divisors at infinity.
For this, we recall that each Segre quartic surface $S\subset\CP_4$
has a structure of a minitwistor space
in the sense of Definition \ref{def:mt}
and therefore the 3-dimensional complex manifold $W_0$
which parameterizes minitwistor lines in $S$
has an Einstein-Weyl structure.
Many Segre surfaces admit double covering structure
over a smooth quadric $Q$ as in Proposition \ref{prop:SS}.
As in Propositions \ref{prop:dcQ} and \ref{prop:si1}.
the dual quadric $Q^*$, considered as included
in the dual space $\CP_4^*$ via the generic projection,
is contained in the completion $W=S^*$ of $W_0$
as a self-intersection locus.
We see that 
{\em the completions of null surfaces in these Einstein-Weyl spaces $W_0$
are always tangent to 
the dual quadrics $Q^*\subset S^*$.}

Before doing so, we first consider the case where
the minitwistor space is a smooth quadric $Q\subset\CP_3$.
In this case, a hyperplane $H\subset\CP_3$ cuts out
a minitwistor line in the original sense iff it is not
tangent to $Q$.
Tangent planes to $Q$ are in one-to-one correspondence
with its tangent point of $Q$,
and thus a hyperplane section $H\subset\CP_3$ does not cut out 
a minitwistor line iff $H\in Q^*$,
where $Q^*\subset\CP_3^*$ is the dual quadric.
This means $W_0^* = \CP_3^*\minus Q^*$.
So we call $Q^*$ a quadric at infinity of $W_0=\CP_3^*\minus Q^*$.
Null surfaces in the Einstein-Weyl space $W_0=\CP_3^*\minus Q^*$
are nothing but the intersections with  the dual planes $p^*\subset\CP_3^*$,
where $p$ is chosen from the minitwistor space $Q$.
For such a dual plane $p^*$, 
the intersection $p^*\cap Q^*$ consists
of the two lines consisting of planes in $\CP_3$ 
which is tangent to $Q$ at some point on the two lines
that pass the point $p$.
This implies that the completion of any null surface in the Einstein-Weyl
space $\CP_3^*\minus Q^*$ is tangent to the infinite quadric $Q^*$.

Returning to the case of genus one,
let $S$ be a Segre quartic surface which has 
a structure of double covering 
over a smooth quadric $Q\subset\CP_3$ as in Proposition \ref{prop:SS}.
As before let $\pi:\CP_4\to\CP_3$ be the generic projection
from the point $w\not\in S$ which induces the double covering map.
By the dual embedding $\pi^*:\CP_3^*\hookrightarrow\CP_4^*$
we identify the dual quadric $Q^*\subset\CP_3^*$ with 
its image $\pi^*(Q^*)$.
We call this a quadric at infinity.
Again null surfaces in $W_0$ are of the form $p^*\cap W_0$,
where $p^*\subset\CP_4^*$ is the dual hyperplane to 
the point $p$ chosen from $S$.
Since the completion $W$ of $W_0$ contains the dual quadric $\pi^*(Q^*)$,
this means that the intersections of the completions of the null surfaces
with the quadric at infinity are of the form
\begin{align}\label{}
p^*\cap \pi^*(Q^*)&= 
\big\{\pi\inv(h)\set h\in Q^*,\,\pi(p)\in h\big\},\quad p\in S.
\end{align}
The RHS can be identified with the union of two lines on $Q\subset\CP_3$
which pass the point $\pi(p)$.
Of course, these two lines intersect transversally at the point $\pi(p)$.
This means that the hyperplane $p^*\subset\CP_4^*$
and the completion of the null surface $p^*\cap W_0$
 is tangent 
to the quadric $\pi^*(Q^*)$ at the point $\pi\inv(T_{\pi(p)} Q)\in \pi^*(Q^*)$.

The second remark is about a relationship between Segre quartic surfaces and
twistor spaces associated
to self-dual metrics on 4-manifolds.
By Jones-Tod \cite{JT85}, a 3-dimensional  Einstein-Weyl manifold is obtained
from a self-dual 4-manifold as a quotient space  under an action of 1-dimensional 
Lie group preserving the self-dual structure.
Correspondingly, a minitwistor space is obtained from the twistor space
of a self-dual 4-manifold as a quotient space with respect to the holomorphic action 
of 1-dimensional complex Lie group.
Among the 16 kinds of the Segre surfaces,
only the ones whose symbol is $[111(11)]$
are known to be obtained this way.
There, the source 4-manifolds are the connected sums of arbitrary number of 
complex projective planes, and their twistor spaces
were constructed in \cite{H2009}.
It might be an interesting question as to whether 
other Segre surfaces can be obtained as quotient spaces of 
compact twistor spaces.

\end{document}